\documentclass[reqno]{amsart}

\usepackage{macros}

\toggletrue{trace}

\makeatletter
\def\enumfix{%
\if@inlabel
 \noindent \par\nobreak\vskip-\topsep\hrule\@height\z@
\fi}

\let\olditemize\itemize
\def\itemize{\enumfix\olditemize}

\makeatother

\begin{document}

\title{The geometry of the cyclotomic trace}

\author{David Ayala, Aaron Mazel-Gee, and Nick Rozenblyum}

\date{\today}

\begin{abstract}
We provide a new construction of the topological cyclic homology $\TC(\cC)$ of any spectrally-enriched $\infty$-category $\cC$, which affords a precise algebro-geometric interpretation of the cyclotomic trace map $\K(X) \ra \TC(X)$ from algebraic K-theory to topological cyclic homology for any scheme $X$.  This construction rests on a new identification of the cyclotomic structure on $\THH(\cC)$, which we find to be a consequence of (i) the geometry of 1-manifolds, and (ii) linearization (in the sense of Goodwillie calculus).  Our construction of the cyclotomic trace likewise arises from the linearization of more primitive data.
\end{abstract}

\maketitle

\setcounter{tocdepth}{2}
\tableofcontents
\setcounter{tocdepth}{2}

\setcounter{section}{-1}

\section{Introduction}
\label{section.intro}

\startcontents[sections]


\subsection{Overview}
\label{subsection.overview}

\introparagraph

\traceismysterious
In this paper, we provide a new construction of $\TC$ which affords a precise conceptual description of the cyclotomic trace at the level of derived algebraic geometry.  The broad idea can be summarized informally as follows.

\begin{slogan}
\label{main.slogan}
The formation of $\TC(X)$ from $\THH(X) \simeq \ms{O}(\ms{L}X)$ amounts to cutting out precisely those functions on the free loopspace $\ms{L}X$
\begin{itemize}
\item which are $\TT$-invariant, and
\item whose value on any loop
\[
(S^1 \xra{\gamma} X)
\in \ms{L}X
\]
determines its values on all its iterates
\[
(S^1 \xra{r} S^1 \xra{\gamma} X) \in \ms{L}X
\]
``up to universal indeterminacy'', subject to all possible universal compatibility relations between iterates of these determinations.
\end{itemize}
Moreover, the cyclotomic trace -- the factorization
\[ \begin{tikzcd}
\K(X)
\arrow{rr}
\arrow[dashed]{rd}
&
&
\THH(X)
\\
&
\TC(X)
\arrow{ru}
\end{tikzcd} \]
of the Dennis trace -- amounts to the recognition that the trace-of-monodromy function of any vector bundle $E \da X$ satisfies these conditions.
\end{slogan}

\noindent The concepts and assertions contained in \Cref{main.slogan} will be both explained and rigorized over the course of this introductory section.

The present paper is the centerpiece of a trilogy, whose overarching purpose is to formulate and prove a precise version of \Cref{main.slogan}.  However, in \Cref{subsection.misc.rmks} we also speculate on some possible applications and future directions of inquiry.  The two supporting papers are \cite{AMR-cyclo}, in which we study cyclotomic spectra (of which $\THH$ is the primary example), and \cite{AMR-fact}, in which we study enriched factorization homology (of which $\THH$ is the primary example).  Our work, particularly our treatment of cyclotomic spectra, finds an important precedent in recent work of Nikolaus--Scholze \cite{NS}; for a discussion of the relationship, see \cite[\Cref*{cyclo:remark.compare.with.NS}]{AMR-cyclo}.

\subsection{The geometry of the Dennis trace}

By definition, the cyclotomic trace arises as an iterated factorization
\begin{equation}
\label{factorizn.of.traces}
\begin{tikzcd}
\K
\arrow{rrr}
\arrow[dashed]{rrd}
\arrow[dashed]{rdd}
&
&
&
\THH
\\
&
&
\THC
\arrow{ru}
& \hspace{-40pt} := \THH^{\htpy \TT}
\\
&
\TC
\arrow{ru}
& \hspace{-35pt} := \THH^{\htpy \Cyclo}
\end{tikzcd}
\end{equation}
of the \textit{Dennis trace} map to \bit{topological Hochschild homology}.  The latter carries a \bit{cyclotomic structure}, i.e.\! it is naturally an object of a certain $\infty$-category $\Cyclo(\Spectra)$ of \textit{cyclotomic spectra}.  This participates in an adjunction
\[ \begin{tikzcd}[column sep=2cm, row sep=0cm]
\Spectra
\arrow[transform canvas={yshift=0.9ex}]{r}{\triv}
\arrow[leftarrow, transform canvas={yshift=-0.9ex}]{r}[yshift=-0.2ex]{\bot}[swap]{(-)^{\htpy \Cyclo}}
&
\Cyclo(\Spectra)
\\
\rotatebox{90}{$\in$}
&
\rotatebox{90}{$\in$}
\\
\TC
&
\THH
\arrow[mapsto]{l}
\end{tikzcd}
~,
\]
whose right adjoint -- which we refer to informally and notationally as the (\textit{homotopy}) \textit{invariants} of the cyclotomic structure -- takes $\THH$ to $\TC$.  Among other data which we will describe shortly, a cyclotomic structure consists of an action of the circle group $\TT$, and we obtain \bit{topological negative cyclic homology} by taking the (homotopy) fixedpoints of the underlying $\TT$-action on $\THH$.

In fact, it is well known how to geometrically interpret the entire upper triangle in diagram \Cref{factorizn.of.traces} \cite{MSV-THH,TV-loops}.  To explain this, suppose that $X$ is a scheme, and let us write $\Perf_X$ for its stable $\infty$-category of vector bundles (i.e.\! perfect complexes).\footnote{In fact, $X$ can be a \textit{spectral} scheme, or more generally any (geometric entity incarnated through its sheaves of some flavor (usually either perfect or coherent) as a) stable $\infty$-category.  In particular, it is of course possible to take $X$ to be an ordinary scheme, but it is only \textit{topological} (i.e.\! $\SS$-linear) Hochschild homology which carries a cyclotomic structure.  (See \cite[\Cref*{cyclo:remark.no.Tate.package.over.Z}]{AMR-cyclo} for more on this point.)}  Then, its algebraic K-theory
\[
\K(X)
:=
\K(\Perf_X)
\]
has cocycles given by vector bundles over $X$, while we have an identification of its topological Hochschild homology
\[
\THH(X)
\simeq
\ms{O}(\ms{L}X)
\]
as the functions on the \bit{free loopspace} of $X$ (a derived mapping stack).\footnote{Actually, the latter identification only holds when $X$ is sufficiently nice (e.g.\! a perfect stack \cite{BZFN}).}  Under these identifications, the Dennis trace is given by the association
\begin{equation}
\label{Dennis.trace}
\begin{tikzcd}[row sep=0cm]
\K(X)
\arrow{r}
&
\THH(X)
\\
\rotatebox{90}{$\in$}
&
\rotatebox{90}{$\in$}
\\
\left(
\begin{array}{c}
\textup{vector bundle}
\\[-0.25cm]
E \da X
\end{array}
\right)
\arrow[mapsto]{r}
&
\left(
\left( \begin{array}{c}
\textup{free loop}
\\[-0.25cm]
S^1 \xra{\gamma} X
\end{array} \right)
\longmapsto
\left(
\begin{array}{c}
\textup{trace of monodromy of}
\\[-0.25cm]
\gamma^* E \da S^1
\end{array}
\right)
\right)
\end{tikzcd}
~.\footnote{After choosing a basepoint of the circle, a vector bundle thereover is determined by its monodromy, an automorphism of the fiber.  Although different basepoints give different automorphisms, these will be conjugate to one another, and hence extracting their trace yields a well-defined invariant.}
\end{equation}
From here, we observe further that the $\TT$-action on $\THH$ corresponds to rotation of loops.  Hence, the factorization of the Dennis trace through $\THC$ in diagram \Cref{factorizn.of.traces} -- which we refer to as the \bit{cyclic trace} -- amounts to recognizing that for any vector bundle $E \da X$, the function obtained by prescription \Cref{Dennis.trace} is invariant under precomposing the loop $\gamma$ with a rotation of the circle.

\subsection{Cocycles for topological Hochschild homology} 
\label{subsection.cocycles}

Let us also take a moment to describe the Dennis trace more directly.

First of all, by definition, $\THH(X)$ is given by the (\bit{enriched}) \bit{factorization homology}
\[
\THH(X)
:=
\THH(\Perf_X)
:=
\int_{S^1} \Perf_X
\]
of the spectrally-enriched $\infty$-category $\Perf_X$ over the framed circle $S^1$.  Thus, informally speaking, a cocycle for $\THH(X)$ is given by the data of
\begin{itemize}
\item a nonempty configuration of points $p_1,\ldots,p_n \in S^1$ on the framed circle (named cyclically),
\item labelings of those points by vectors bundles $E_i \da X$, and
\item labelings of the framed intervals between them by maps $E_i \ra E_{i+1}$.
\end{itemize}
The equivalence relation among cocycles comes from allowing points
\begin{itemize}
\item to disappear, in which case we compose the adjacent maps, and
\item to anticollide, in which case we label the new interval by the identity map.
\end{itemize}

Somewhat more precisely, $\THH(X)$ is given by the colimit of a diagram indexed over such configurations and point-labelings, with value the smash product spectrum
\begin{equation}
\label{value.of.diagram.for.THH}
\bigotimes_{i \in \ZZ/n} \ulhom_{\Perf_X}(E_i,E_{i+1})
~.
\end{equation}
(Of course the structure maps of this diagram are just as described above, only they use \textit{enriched} composition and unit maps of $\Perf_X$.)

Note that among these two heuristic descriptions of $\THH(X)$, the former would suggest that we are simply taking a cartesian product of hom-spectra in $\Perf_X$;\footnote{In fact, this is really more of a description of \textit{unenriched} factorization homology, which in place of the smash product \Cref{value.of.diagram.for.THH} of hom-spectra uses cartesian products of hom-spaces.} by contrast, the smash product \Cref{value.of.diagram.for.THH} enforces \textit{multilinearity}.  This distinction will be crucial.

At the level of cocycles, the Dennis trace takes a vector bundle $E \da X$ to the cocycle for $\THH(X)$ given by a single point on the circle labeled by $E$, with the framed interval labeled by its identity map.  The identification $\THH(\Perf_X) \simeq \ms{O}(\ms{L}X)$ passes through \textit{Morita invariance}, which yields an equivalence
\[
\THH(\Perf_A)
\simeq
\THH(\fB A)
\]
for any ring spectrum $A$, where $\fB A$ denotes its delooping to a one-object spectrally-enriched $\infty$-category (see Remarks \ref{flagged} \and \ref{morita.invce}).

\subsection{The geometry of the cyclotomic trace}
\label{subsection.geometry}

Whereas the cyclic trace arises from \textit{automorphisms} of the circle, the cyclotomic trace arises from \textit{endomorphisms} of the circle -- that is, from the map
\[ \begin{tikzcd}[row sep=0cm]
\ms{L}X
\arrow{r}{r^*}
&
\ms{L}X
\\
\rotatebox{90}{$\in$}
&
\rotatebox{90}{$\in$}
\\
\left( S^1 \xra{\gamma} X \right)
\arrow[mapsto]{r}
&
\left( S^1 \xra{r} S^1 \xra{\gamma} X \right)
\end{tikzcd} \]
given by precomposing with an $r$-fold self-covering map, for a natural number $r \in \Nx$.  Namely, given a vector bundle $E \da X$, there is a sense in which the value of its trace-of-monodromy function on a loop $\gamma$ determines its value on the loop $r^* \gamma$: not quite on-the-nose, but \textit{up to universal indeterminacy}.

In fact, this phenomenon already manifests itself in the case of a square matrix $M \in {\sf M}_{d \times d}(R)$ over a ring $R$; if this matrix records the monodromy around $\gamma$, then its $r\th$ power records the monodromy around $r^*\gamma$.  Let us therefore illustrate this phenomenon by considering the two elements
\begin{equation}
\label{two.elts.of.rth.tensor.power}
\tr \left( M^{\otimes r} \right)
,
\tr(M)^{\otimes r}
\in
R^{\otimes r}
\end{equation}
of the tensor power of the base ring.  We first observe that the elements \Cref{two.elts.of.rth.tensor.power} are both cyclically invariant, the latter manifestly and the former by standard properties of traces of matrices; in other words, they actually lie in the fixedpoints
\[ (R^{\otimes r})^{\Cyclic_r} \subset R^{\otimes r} \]
with respect to the action of the cyclic group $\Cyclic_r$.  Now, the elements \Cref{two.elts.of.rth.tensor.power} are not equal, but they differ by \textit{norms}.  The first nontrivial example occurs when $r=2$.  In this case, we have the \bit{norm map}
\[ \begin{tikzcd}[row sep=0cm]
(R^{\otimes 2})_{\Cyclic_2}
\arrow{r}{\Nm}
&
(R^{\otimes 2})^{\Cyclic_2}
\\
\rotatebox{90}{$\in$}
&
\rotatebox{90}{$\in$}
\\
{[x \otimes y]}
\arrow[mapsto]{r}
&
\sum_{\sigma \in \Cyclic_2} \sigma(x \otimes y)
\end{tikzcd} \]
from the orbits to the fixedpoints, and we can compute directly that the difference
\begin{align*}
\tr \left( M^{\otimes 2} \right) - \tr(M)^{\otimes 2}
&=
\sum_{i,j} m_{i,j} \otimes m_{j,i}
-
\sum_{k,l} m_{k,k} \otimes m_{l,l}
\\
&=
\sum_{i \not= j} m_{i,j} \otimes m_{j,i}
-
\sum_{k \not= l} m_{k,k} \otimes m_{l,l}
\\
&=
\sum_{i<j} \Nm ( [ m_{i,j} \otimes m_{j,i} ] )
-
\sum_{k<l} \Nm ( [ m_{k,k} \otimes m_{l,l} ] )
\end{align*}
does indeed land in its image.  Thus, we obtain a canonical identification
\[ \tr(M^{\otimes 2}) \equiv \tr(M)^{\otimes 2} \]
in the \textit{cokernel} of the norm map.  More generally, the elements \Cref{two.elts.of.rth.tensor.power} differ by norms not from just the trivial subgroup of $\Cyclic_r$ but from \textit{all} proper subgroups; this assertion categorifies the combinatorial observation that for instance if
\[ M = \left( \begin{array}{cc} m_{1,1} & 0 \\ 0 & m_{2,2} \end{array} \right) \]
is a $2 \times 2$ diagonal matrix, then the difference $\tr(M^{\otimes r}) - \tr(M)^{\otimes r}$ is governed by the binomial coefficients ${r \choose i}$ for $0 < i < r$, which are never coprime to $r$.\footnote{For a more thorough and coordinate-independent treatment of these phenomena, we refer the reader to \cite{NS}.}

Let us now return to considering $\THH(X)$, i.e.\! functions on the free loopspace of $X$.  In this setting, we would like to compare two elements of
\[
\THH(X)^{\tate \Cyclic_r}
~,
\]
the (``\bit{generalized}'', though we'll omit this for brevity) \bit{Tate construction} for the natural $\Cyclic_r$-action on $\THH(X)$ -- a homotopical version of the ``quotient the fixedpoints by norms from all proper subgroups'' operation.\footnote{When $r=p$ is prime, then the only proper subgroup of $\Cyclic_p$ is the trivial subgroup.  Correspondingly, in this case the generalized Tate construction $(-)^{\tate \Cyclic_p}$ reduces to the ordinary Tate construction, namely the cofiber
\[
(-)_{\htpy \Cyclic_p}
\xra{\Nm}
(-)^{\htpy \Cyclic_p}
\longra
(-)^{{\sf t} \Cyclic_p}
~.
\]
}\footnote{For a justification of this description of the generalized Tate construction, see \cite[\Cref*{cyclo:remark.genzd.tate}]{AMR-cyclo}.}  As we will soon see, this indeed receives two distinguished maps: one running
\begin{equation}
\label{intro.can.str.map}
\THH(X)^{\htpy \Cyclic_r}
\longra
\THH(X)^{\tate \Cyclic_r}
\end{equation}
and corresponding to the construction ``$M \mapsto \tr(M)^{\otimes r}$'', and another running
\begin{equation}
\label{intro.cyclo.str.map}
\THH(X)
\longra
\THH(X)^{\tate \Cyclic_r}
\end{equation}
and corresponding to the construction ``$M \mapsto \tr \left( M^{\otimes r} \right)$''.\footnote{It follows from an analysis at the level of cyclic bar constructions that these maps do indeed behave as claimed.  The key observation is that the trace actually occurs in the equivalence $\THH( \fB({\sf M}_{d \times d}(R))) \xra{\sim} \THH(\fB R)$ coming from Morita invariance (see e.g.\! \cite[Theorem 1.2.4]{Loday-cyclic}); then, we can see the map \Cref{intro.can.str.map} as coming from the composite
\[
\left( {\sf M}_{d \times d}(R)^{\otimes r(\bullet+1)} \right)^{\htpy \Cyclic_r}
\longra
\left( R^{\otimes r(\bullet+1)} \right)^{\htpy \Cyclic_r}
\longra
\left( R^{\otimes r(\bullet+1)} \right)^{\tate \Cyclic_r}
\]
and the map \Cref{intro.cyclo.str.map} as coming from the composite
\[
{\sf M}_{d \times d}(R)^{\otimes (\bullet+1)}
\longra
\left( {\sf M}_{d \times d}(R)^{\otimes r (\bullet+1)} \right)^{\tate \Cyclic_r}
\longra
\left( R^{\otimes r(\bullet+1)} \right)^{\tate \Cyclic_r}
~.
\]
(More generally we might consider arbitrary composable (not necessarily square) matrices, but this would only serve to make the notation more complicated, and is itself rendered unnecessary by Morita invariance; the prior reduction to free objects (instead of just their retracts) is likewise justified by Morita invariance.)
}  These are suitably equivariant, in such a way that they determine a pair of parallel maps
\[
\begin{tikzcd}
\THH(X)^{\htpy \TT}
\arrow[transform canvas={yshift=0.8ex}]{r}
\arrow[transform canvas={yshift=-0.8ex}]{r}
&
\left( \THH(X)^{\tate \Cyclic_r} \right)^{\htpy \TT}
\end{tikzcd}
~;
\]
unwinding the definitions, it follows that the equalizer of this pair picks out
\begin{itemize}
\item[] \textit{$\TT$-equivariant functions on $\ms{L}X$ whose value on any loop $\gamma$ determines its values on the iterate $r^*\gamma$ up to universal indeterminacy},
\end{itemize}
thus embodying much of the algebro-geometric description of $\TC(X)$ asserted in \Cref{main.slogan}.

But much more is true: these structure maps will come equipped with an \textit{infinite hierarchy} of compatibilities, corresponding to all possible composable sequences of self-coverings of the circle as well as a choice at each stage of whether to enact the operation ``$M \mapsto \tr(M)^{\otimes r}$'' or the operation ``$M \mapsto \tr \left(M^{\otimes r} \right)$''.\footnote{In fact, there will be yet more compatibilities, accounting for the fact that the ``quotienting by norms from all proper subgroups'' operation doesn't strictly commute with multiplication of natural numbers (the orders of the various cyclic groups); see \Cref{tau.action.only.left.lax}.}  Collectively, the limit over all of these structure maps will impose
\begin{itemize}
\item[] \textit{all possible universal compatibility relations between iterates of these determinations},
\end{itemize}
which is the remaining assertion of \Cref{main.slogan} regarding the geometry of $\TC(X)$; these compatibilities will be spelled out more explicitly in \Cref{geometry.of.lim.sd.BN}.  All in all, we will see that in an admittedly revisionistic sense,
\begin{itemize}
\item[] \textit{$\TC(X)$ is precisely built to encode all of the structure present on those functions on $\ms{L}X$ that arise as trace-of-monodromy functions of vector bundles over $X$},
\end{itemize}
completing the explanation of \Cref{main.slogan}.

\subsection{Main theorems}

We now describe the four central theorems of the trilogy to which this paper belongs, which collectively implement the geometric picture described in \Cref{subsection.geometry} above; the first and third are proved in \cite{AMR-cyclo}, while the second and fourth are proved here.

First of all, we reidentify the $\infty$-category $\Cyclo(\Spectra)$ of cyclotomic spectra in terms of the Tate construction.

\begin{maintheorem}[A reidentification of cyclotomic spectra (proved in \cite{AMR-cyclo})]
\label{cyclo.is.rlax.lim}
The endofunctors
\[ \begin{tikzcd}[column sep=2cm, row sep=0cm]
\Fun(\BT,\Spectra)
\arrow{r}
&
\Fun({\sf B}(\TT/\Cyclic_r),\Spectra)
\arrow{r}{(\TT/\Cyclic_r) \simeq \TT}[swap]{\sim}
&
\Fun(\BT,\Spectra)
\\
\rotatebox{90}{$\in$}
&
\rotatebox{90}{$\in$}
&
\rotatebox{90}{$\in$}
\\
(\TT \lacts E)
\arrow[mapsto]{r}
&
\left(\TT/\Cyclic_r \lacts E^{\tate \Cyclic_r}\right)
\arrow[mapsto]{r}
&
\left(\TT \lacts E^{\tate \Cyclic_r}\right)
\end{tikzcd} \]
for $r \in \Nx$ assemble into a left-lax right action
\[
\Fun(\BT,\Spectra)
\ractstau
\Nx
~.
\]
Moreover, there is a canonical equivalence of $\infty$-categories
\[
\Cyclo(\Spectra)
\simeq
\lim^\rlax \left( \Fun(\BT,\Spectra) \ractstau \Nx \right)
\]
with the right-lax limit of this left-lax right action.
\qed
\end{maintheorem}

\begin{remark}
\label{tau.action.only.left.lax}
To say that the action is \textit{left-lax} is to say that there are only comparison morphisms (rather than equivalences) governing its associativity.  For instance, for each pair of elements $r,s \in \Nx$ and any $E \in \Fun(\BT,\Spectra)$, there is a canonical map
\begin{equation}
\label{comparison.map.chi}
E^{\tate \Cyclic_{rs}}
\xra{\chi_{r,s}}
\left( E^{\tate \Cyclic_r} \right)^{\tate \Cyclic_s}
~.\footnote{More generally, for any word $W := (r_1,\ldots,r_n)$ of elements of $\Nx$, writing $|W| := \prod_i r_i$ for the product of its elements, there is a canonical map
\[
E^{\tate \Cyclic_{|W|}}
\xra{\chi_W}
E^{\tate \Cyclic_W}
:=
E^{\tate \Cyclic_{r_1} \cdots \tate \Cyclic_{r_n}}
~.
\]
These maps satisfy various compatibilities, a few of which will appear presently.}
\end{equation}
Then, an object of the \textit{right-lax} limit of this left-lax right action -- that is, a cyclotomic spectrum -- consists of the following data:
\begin{enumerate}

\setcounter{enumi}{-1}

\item an object
\[
T
\in
\Fun(\BT,\Spectra)
~,
\]
namely the \bit{underlying $\TT$-spectrum} of the cyclotomic spectrum;

\item for every element $r \in \Nx$, a \bit{cyclotomic structure map}
\[
T
\xlongra{\sigma_r}
T^{\tate \Cyclic_r}
\]
in $\Fun(\BT,\Spectra)$;

\item for every pair of elements $r,s \in \Nx$, the \textit{data} of a commutative square
\[ \begin{tikzcd}[row sep=1.5cm, column sep=1.5cm]
T
\arrow{r}{\sigma_s}
\arrow{d}[swap]{\sigma_{rs}}
&
T^{\tate \Cyclic_s}
\arrow{d}{(\sigma_r)^{\tate \Cyclic_s}}
\\
T^{\tate \Cyclic_{rs}}
\arrow{r}[swap]{\chi_{r,s}}
&
\left( T^{\tate \Cyclic_r} \right)^{\tate \Cyclic_s}
\end{tikzcd} \]
in $\Fun(\BT,\Spectra)$, denoted $\sigma_{r,s}$;

\item for every triple of elements $r,s,t \in \Nx$, the \textit{data} of a commutative 3-cube
\[ \begin{tikzcd}[row sep=1.5cm]
&
T^{\tate \Cyclic_{rst}}
\arrow{rr}{\chi_{rs,t}}
\arrow{dd}[pos=0.3]{\chi_{r,st}}
&
&
\left( T^{\tate \Cyclic_{rs}}\right)^{\tate \Cyclic_{t}}
\arrow{dd}{\left( \chi_{r,s} \right)^{\tate \Cyclic_{t}}}
\\
T
\arrow[crossing over]{rr}[pos=0.7]{\sigma_{t}}
\arrow{ru}{\sigma_{rst}}
\arrow{dd}[swap]{\sigma_{st}}
&
&
T^{\tate \Cyclic_{t}}
\arrow{ru}[swap]{\left( \sigma_{rs} \right)^{\tate \Cyclic_{t}}}
\\
&
\left( T^{\tate \Cyclic_{r}} \right)^{\tate \Cyclic_{st}}
\arrow{rr}[pos=0.3]{\chi_{s,t}}
&
&
\left( \left( T^{\tate \Cyclic_{r}} \right)^{\tate \Cyclic_{s}} \right)^{\tate \Cyclic_{t}}
\\
T^{\tate \Cyclic_{st}}
\arrow{rr}[swap]{\chi_{s,t}}
\arrow{ru}[swap]{\left( \sigma_{r} \right)^{\tate \Cyclic_{st}}}
&
&
\left( T^{\tate \Cyclic_{s}} \right)^{\tate \Cyclic_{t}}
\arrow{ru}[swap]{\left( \left( \sigma_{r} \right)^{\tate \Cyclic_{s}} \right)^{\tate \Cyclic_{t}}}
\arrow[leftarrow, crossing over]{uu}[pos=0.7, swap]{\left( \sigma_{s} \right)^{\tate \Cyclic_{t}}}
\end{tikzcd} \]
in $\Fun(\BT,\Spectra)$, denoted $\sigma_{r,s,t}$,
\begin{itemize}
\item whose front, top, and left, and right faces are the commutative squares $\sigma_{s,t}$, $\sigma_{rs,t}$, $\sigma_{r,st}$, and $(\sigma_{r,s})^{\tate \Cyclic_t}$, respectively,
\item whose bottom face is determined by the naturality of $\chi_{s,t}$, and 
\item whose back face is determined by the compatibility of the various $\chi$ maps \Cref{comparison.map.chi};
\[ \vdots \]
\end{itemize}
\end{enumerate}

\begin{enumerate}[label={($\alph*$)}]

\setcounter{enumi}{13}

\item in general, for each word $W := ( r_1,\ldots,r_n) $ in $\Nx$, the \textit{data} of a commutative $n$-cube $\sigma_W$ in $\Fun(\BT,\Spectra)$, whose faces are given
\begin{itemize}
\item by previously specified data (inductively for all $n \geq 0$) and
\item by the various coherences enjoyed by the ``$\chi$'' comparison maps \Cref{comparison.map.chi}.
\end{itemize}

\end{enumerate}
These data are subject to certain compatibility relations that apply whenever the unit element $1 \in \Nx$ appears in a word.
\end{remark}

\begin{remark}
The sorts of words in $\Nx$ that are relevant to us are parametrized by the subdivision category $\sd(\BN)$: its objects are precisely the equivalence classes of words in $\Nx$ under the relation that any instances of the element $1 \in \Nx$ may be freely inserted or omitted.\footnote{Due to the equivalence relation, it doesn't matter mathematically whether we allow the empty word $\es$ or not, since it's equivalent to the word $(1)$.  However, it will be convenient to use the evident ``normal form'' where we simply omit all $1$'s, which leads us to choose to allow the empty word.}  In fact, this is not a coincidence: subdivision plays a key role in the definition of a \textit{right}-lax limit of a \textit{left}-lax action (i.e.\! it appears because these handednesses disagree), and it will reappear in our explicit formula for $\TC$ (see \Cref{TC.formula} and \Cref{geometry.of.lim.sd.BN}).
\end{remark}

In order to state our next theorem, we introduce/recall the following notation:
\begin{itemize}
\item $\Spaces$ and $\Spectra$ respectively denote the $\infty$-categories of spaces and of spectra;
\item $\Cat(\Spaces) \simeq \Cat$ and $\Cat(\Spectra)$ respectively denote the $\infty$-categories of spatially-enriched (i.e.\! unenriched) $\infty$-categories and of spectrally-enriched $\infty$-categories;
\item $\THH_\Spaces$ and $\THH:= \THH_\Spectra$ respectively denote spatially-enriched and spectrally-enriched factorization homology over the circle;
\item $\Cyclo^\htpy(\Spaces)$ and $\Cyclo(\Spectra) := \Cyclo^\tate(\Spectra)$ respectively denote the $\infty$-categories of unstable cyclotomic spaces (recalled in \Cref{section.roadmap.unst.cyclo.str} below) and of cyclotomic spectra.
\end{itemize}
 
\needspace{2\baselineskip}
\begin{maintheorem}[The cyclotomic structure on $\THH$ (proved in \Cref{part.cyclo.str})]
\label{THH.is.cyclo}
\begin{enumerate}
\item[]
\item There is a \bit{diagonal package} for spaces, which induces a canonical lift
\[ \begin{tikzcd}[column sep=1.5cm, row sep=1.5cm]
\Cat(\Spaces)
\arrow[dashed]{r}
\arrow{rd}[swap]{\THH_\Spaces}
&
\Cyclo^\htpy(\Spaces)
\arrow{d}{\fgt}
\\
&
\Spaces
\end{tikzcd}
\]
through the forgetful functor, yielding the \bit{unstable cyclotomic structure} on $\THH_\Spaces$.
\item Via linearization (in the sense of Goodwillie calculus), the diagonal package for spaces gives rise to the \bit{Tate package} for spectra, which induces a canonical lift
\[ \begin{tikzcd}[column sep=1.5cm, row sep=1.5cm]
\Cat(\Spectra)
\arrow[dashed]{r}
\arrow{rd}[swap]{\THH}
&
\Cyclo(\Spectra)
\arrow{d}{\fgt}
\\
&
\Spectra
\end{tikzcd}
\]
through the forgetful functor, yielding the \bit{cyclotomic structure} on $\THH$.
\end{enumerate}
\end{maintheorem}

\begin{remark}
The term \bit{package} is meant to emphasize that this object \textit{packages} together the named maps, along with all of their requisite functorialities: the \textit{diagonal package} packages the diagonal maps
\[
V
\longra
\left( V^{\times r} \right)^{\htpy \Cyclic_r}
\]
for spaces $V \in \Spaces$ (which are in fact equivalences), while the \textit{Tate package} packages the \bit{Tate diagonal} maps
\[
E
\longra
\left( E^{\otimes r} \right)^{\tate \Cyclic_r}
\]
for spectra $E \in \Spectra$.  The way in which diagonal maps for spaces induce the unstable cyclotomic structure on $\THH_\Spaces$ is described in \Cref{subsection.unwind.unst.cyclo.str.map}; the Tate diagonal maps for spectra induce the cyclotomic structure on $\THH$ in an essentially identical way.
\end{remark}

Together, \Cref{cyclo.is.rlax.lim,THH.is.cyclo} go much of the way towards formalizing the assertions of \Cref{main.slogan} regarding $\TC$.  However, what still lacks is a concrete formula for the passage from $\THH$ to $\TC$, guaranteeing that it does indeed implement the asserted geometric description.  Our next result fills this gap.

\begin{maintheorem}[The formula for $\TC$ (proved in \cite{AMR-cyclo})]
\label{TC.formula}
There is a canonical factorization
\[ 
\begin{tikzcd}[row sep=1.5cm]
\Spectra
\arrow[leftarrow]{rr}{(-)^{\htpy \Cyclo}}
&
&
\Cyclo(\Spectra)
\arrow{ld}{(-)^{\htpy \TT}}
\\
&
\Fun(\sd(\BN),\Spectra)
\arrow{lu}{\lim}
\end{tikzcd}
~,\footnote{It would be slightly more correct to write $\sd \left( \BNop \right)$ instead of $\sd(\BN)$, but these two categories are canonically equivalent (due to the commutativity of the monoid $\Nx$) and so to simplify our notation we elide the distinction here.}
\]
where the first functor takes a cyclotomic spectrum $T \in \Cyclo(\Spectra)$ to the diagram
\begin{equation}
\label{diagram.indexed.by.sd.BN}
\begin{tikzcd}[row sep=0cm, column sep=1.5cm]
&
\sd(\BN)
\arrow{r}{T^{\htpy \TT}}
&
\Spectra
\\
&
\rotatebox{90}{$\in$}
&
\rotatebox{90}{$\in$}
\\
(r_1,\ldots,r_k)
=:
\hspace{-2.2cm}
&
W
\arrow[maps to]{r}
&
\left( T^{\tate \Cyclic_W} \right)^{\htpy \TT}
&
\hspace{-1.7cm}
:= \left( T^{\tate \Cyclic_{r_1} \cdots \tate \Cyclic_{r_k}} \right)^{\htpy \TT}
\end{tikzcd}
\end{equation}
of spectra.
\qed
\end{maintheorem}

\begin{remark}
\label{geometry.of.lim.sd.BN}
The category $\sd(\BN)$ and diagram \Cref{diagram.indexed.by.sd.BN} are described in detail in \cite[\Cref*{cyclo:unwind.sd.BN}]{AMR-cyclo}.  As described there, the morphisms in $\sd(\BN)$ are generated by the operations of
\begin{itemize}
\item adding a letter to the beginning of a word,
\item adding a letter to the end of a word, and
\item factoring a letter in a word.
\end{itemize}
In particular, the empty word
\[
\emptyword \in \sd(\BN)
\]
is an empty object (i.e.\! it receives no nonidentity maps), and moreover it maps to every other object in at least two distinct ways.  Hence, the extraction of the limit
\[
\lim \Cref{diagram.indexed.by.sd.BN}
=:
T^{\htpy \Cyclo}
\]
amounts to \textit{cutting out equations} inside of the value
\[
T^{\htpy \TT}
\in
\Spectra
\]
of the diagram \Cref{diagram.indexed.by.sd.BN} at the empty word.

Let us now specialize to our particular case of interest, namely when we take our cyclotomic spectrum to be
$
T
=
\THH(X)
$
for a scheme $X$.  In this case, we have that
\[
\THH(X)^{\htpy \TT}
\simeq
\left( \ms{O}(\ms{L}X) \right)^{\htpy \TT}
\]
consists of the $\TT$-invariant functions on the free loopspace $\ms{L}X$.  Then, the three operations described above can be interpreted \textit{geometrically}: they correspond to operations on such functions of the form
\begin{itemize}
\item ``$M \mapsto \tr \left( M^{\otimes r} \right)$'',
\item ``$M \mapsto \tr(M)^{\otimes r}$'', and
\item ``take a further quotient corresponding to the new universal indeterminacies that arise in factoring an $rs$-fold covering map as an $r$-fold covering map followed by an $s$-fold covering map''.
\end{itemize}
Thus, \Cref{TC.formula} shows that the passage from $\THH(X)$ to $\TC(X)$ afforded by \Cref{THH.is.cyclo} does indeed implement the geometric description asserted in \Cref{main.slogan}.
\end{remark}

In order to state our final theorem, we introduce/recall the following notation.
\begin{itemize}
\item We write
\[
\TC_\Spaces
:
\Cat
\xra{\THH_\Spaces}
\Cyclo^\htpy(\Spaces)
\xra{(-)^{\htpy\WW}}
\Spaces
\]
for the unstable topological cyclic homology functor (recalled in \Cref{section.unstable.cyclo.trace}).  This is the composite of spatially-enriched factorization homology over the circle with the ``homotopy invariants of the unstable cyclotomic structure'' functor; the notation $\WW$ is explained in \Cref{section.roadmap.unst.cyclo.str}.
\item We write
\[
\TC
:
\Cat(\Spectra)
\xra{\THH}
\Cyclo(\Spectra)
\xra{(-)^{\htpy \Cyclo}}
\Spectra
\]
for the topological cyclic homology functor.  This is the composite of spectrally-enriched factorization homology over the circle with the ``homotopy invariants of the cyclotomic structure'' functor.
\item For any enriching $\infty$-category $\cV$, we write
\[
\Cat(\cV)
\xlongra{\iota}
\Spaces
\]
for the ``underlying $\infty$-groupoid'' functor.
\end{itemize}

\needspace{2\baselineskip}
\begin{maintheorem}[The cyclotomic trace (proved in \Cref{part.cyclo.trace})]
\label{cyclo.trace}
\begin{enumerate}
\item[]
\item The contravariant functoriality of spatially-enriched factorization homology with respect to the proper constructible bundle
\[
\DD^0
\longla
S^1
\]
induces an \bit{unstable Dennis trace} 
\[
\iota
\simeq \int_{\DD^0}
\longra
\int_{S^1}
=:
\THH_\Spaces
\]
in $\Fun(\Cat,\Spaces)$, whose equivariance induces a factorization
\[ \begin{tikzcd}
\iota
\arrow{rr}
\arrow[dashed]{rd}
&
&
\THH_\Spaces
\\
&
\TC_\Spaces
\arrow{ru}
\end{tikzcd} \]
through the \bit{unstable cyclotomic trace}.

\item Via linearization (in the sense of Goodwillie calculus) and the equivalence
\[
\Sigma^\infty_+ \iota
\simeq
\int_{\DD^0}
\]
in $\Fun(\Cat(\Spectra),\Spectra)$, the unstable Dennis and cyclotomic trace maps give rise to the \bit{Dennis pre-trace} and \bit{cyclotomic pre-trace} maps
\[ \begin{tikzcd}[row sep=1.5cm, column sep=1.5cm]
\Sigma^\infty_+ \iota
\arrow[dashed]{r}
\arrow[dashed]{rd}
&
\THH
\\
&
\TC
\arrow{u}
\end{tikzcd} \]
in $\Fun(\Cat(\Spectra),\Spectra)$, the latter of which when restricted to the subcategory
\[
\StCat
\subset
\Cat(\Spectra)
\]
of stable $\infty$-categories and exact functors between them induces the \bit{cyclotomic trace} map
\[ \begin{tikzcd}[row sep=1.5cm, column sep=1.5cm]
\Sigma^\infty_+ \iota
\arrow{r}
\arrow{rd}
\arrow{d}
&
\THH
\\
\K
\arrow[dashed]{r}
&
\TC
\arrow{u}
\end{tikzcd} \]
in $\Fun(\StCat,\Spectra)$.

\end{enumerate}
\end{maintheorem}

\begin{remark}
It is immediate that for a scheme $X$, the Dennis pre-trace
\[
\Sigma^\infty_+ \iota (\Perf_X)
\longra
\THH(\Perf_X)
=:
\THH(X)
\]
behaves as described in \Cref{subsection.cocycles}: it takes a vector bundle $E \in \Perf_X$ to the cocycle for $\THH(X)$ represented by a single point on the circle labeled by $E$, with the framed interval labeled by its identity map.  As under the identification
\[
\THH(X)
\simeq
\ms{O}(\ms{L}X)
\]
this cocycle corresponds to the trace-of-monodromy function of $E$, the cyclotomic (pre-)trace indeed amounts to the observation that such a function necessarily satisfies the conditions imposed by taking the homotopy invariants of the cyclotomic structure.  Thus, \Cref{cyclo.trace} completes the rigorization of the assertions of \Cref{main.slogan}.
\end{remark}

\subsection{Miscellaneous remarks}
\label{subsection.misc.rmks}

\begin{remark}
There is an evident cardinality filtration
\[
\THH
\simeq
\colim
\left(
\THH_{\leq 1}
\longra
\THH_{\leq 2}
\longra
\THH_{\leq 3}
\longra
\cdots
\right)
\]
of $\THH$, obtained by restricting the number of points on the circle.  This filtration is respected by the $\TT$-action, and so it refines to a filtration
\[
(\THC)_{\leq n}
:=
\left( \THH_{\leq n} \right)^{\htpy \TT}
\]
of
$\THC := \THH^{\htpy \TT}$.\footnote{This will not generally be ``exhaustive'' (i.e.\! it won't generally recover $\THC$ in the colimit), since directed colimits don't generally commute with homotopy $\TT$-fixedpoints.}  In fact, our construction provides for a further refinement to a filtration of $\TC$: for a word
\[
W
:=
(r_1,\ldots,r_k)
\in
\sd(\BN)
~,
\]
if we write
\[
|W|
:=
\prod_i r_i
\in
\Nx
\]
for the product of its elements, then we set
\[
\TC_{\leq n}
:=
\lim_{\sd(\BN)}
\left( (\THH_{\leq n \cdot |W|})^{\tate \Cyclic_W} \right)^{\htpy \TT}
~.\footnote{Note that this filtration of $\TC$ does not arise from a filtration of $\THH$ by cyclotomic spectra.}
\]
Preliminary computations suggest that these filtrations should be closely related to the Goodwillie towers of the \textit{relative} versions of the respective theories.  This would give purchase on nonconnective and nonaffine generalizations of the Dundas--Goodwillie--McCarthy theorem, which at present only applies in the somewhat restrictive case of connective ring spectra.\end{remark}

\begin{remark}
Whereas most previous work on $\THH$ and $\TC$ relies on simplicial methods (e.g.\! the cyclic bar construction), our work uses the geometry of 1-manifolds.  This brings conceptual advantages: it clarifies the connection with loopspaces, and it suggests the appropriate higher-dimensional generalization (see \Cref{higher.TC}).  But it also brings an important technical advantage: the desired symmetries are built into the framework, rather than having to be put back in by hand (as e.g.\! an action of the simplicial circle).  In effect, employing simplicial methods amounts to choosing a basepoint on the circle (see \cite[\Cref*{fact:excision.for.enr.fact.hlgy}]{AMR-fact}), which destroys much of its symmetry.
\end{remark}

\begin{remark}
\label{higher.TC}
Our work paves the way for \textit{higher-dimensional} versions of $\TC$, defined for any spectrally-enriched $(\infty,d)$-category: this should simultaneously involve all compact framed $d$-manifolds and (framed) covering maps among them.  One of our main technical results -- the \textit{proto Tate package} \cite[\Cref*{cyclo:proto.tate}]{AMR-cyclo}, which mediates the passage from the diagonal package to the Tate package of \Cref{THH.is.cyclo} -- is already sufficiently general to apply in this situation, since it handles arbitrary finite groups and their iterated extensions (appearing as the deck groups of covering maps).  We expect that this should reduce to the notion of ``covering homology'' \cite{BCD-cov} in the case of a connective commutative ring spectrum (viewed as a flagged (see \Cref{flagged}) spectrally-enriched $(\infty,d)$-category with a single $k$-morphism for all $k<d$).
\end{remark}

\begin{remark}
We phrase our work in terms of $\infty$-categories enriched in spaces and spectra, but in fact it applies word-for-word to the more general context of a (presentably) \textit{cartesian} symmetric monoidal $\infty$-category $\cV$ and its stabilization ${\sf Stab}(\cV)$.  On the other hand, much of our development of enriched factorization homology is stated in terms of an \textit{arbitrary} symmetric monoidal $\infty$-category $\cV$; we only specialize to spaces and spectra when we start to use their more particular features.
\end{remark}

\begin{remark}
\label{flagged}
We work not just with $\cV$-enriched $\infty$-categories, but with \textit{flagged} $\cV$-enriched $\infty$-categories, i.e.\! ``a $\cV$-enriched $\infty$-category with a chosen space of objects''.\footnote{\label{complete.V.enr.cats}An enriched $\infty$-category $\cC$ has a canonically defined maximal subgroupoid, but one can also present it as having for its  ``space of objects'' some other $\infty$-groupoid: to do so amounts to choosing a $\pi_0$-surjection from the latter to the former.  Of course, taking this map to be an equivalence recovers the terminal such choice: this generalizes the relationship between Segal spaces and \textit{complete} Segal spaces (see \cite[\Cref*{fact:remark.segal.spaces}]{AMR-fact}).}  In particular, a ring spectrum $A$ determines a flagged spectrally-enriched $\infty$-category $\fB A$ whose ``space of objects'' consists of a single point.\footnote{Indeed, these two types of data are manifestly equivalent; see \Cref{subsection.enr.infty.cats}.}  This added generality is useful (and requires no extra work), since it recovers the more classical definition of $\THH(A)$ via the cyclic bar construction (see \cite[\Cref*{fact:excision.for.enr.fact.hlgy}]{AMR-fact}).
\end{remark}

\begin{remark}
\label{morita.invce}
Enriched factorization homology over the circle satisfies Morita invariance.  To illustrate this, we offer the following proof-by-picture.  Let $\cC$ and $\cD$ be flagged $\cV$-enriched $\infty$-categories, and suppose that
\[
\cC^\op \times \cD
\xlongra{\cP}
\cV
\]
and
\[
\cD^\op \times \cC
\xlongra{\cQ}
\cV
\]
are bimodules determining a Morita equivalence between $\cC$ and $\cD$.  Then, we have the string of equivalences
\[
\begin{tikzpicture}

\draw (-0.5,0) node {\scalebox{3}{$\int$}};
\draw[
        decoration={markings, mark=at position 0.5 with {\arrow{>}}},
        postaction={decorate}
        ] (0.5,0) circle (0.5cm);
\draw (0,0) node[anchor=east] {$\cC$};

\draw (2,0) node {$\simeq$};

\draw (3,0) node {\scalebox{3}{$\int$}};
\draw[
        decoration={markings, mark=at position 0.5 with {\arrow{>}}},
        postaction={decorate}
        ] (4,0) circle (0.5cm);
\draw (3.5,0) node[anchor=east] {$\cC$};
\draw (4.5,0) node {$\bullet$};
\draw (4.5,0) node[anchor=west] {$\cC$};

\draw (5.5,0) node {$\simeq$};

\draw (6.5,0) node {\scalebox{3}{$\int$}};
\draw[
        decoration={markings, mark=at position 0.5 with {\arrow{>}}},
        postaction={decorate}
        ] (7.5,0) circle (0.5cm);
\draw (7,0) node[anchor=east] {$\cC$};
\draw (8,0) node {$\bullet$};
\draw (8,0) node[anchor=west] {$(\cQ \otimes_\cD \! \cP)$};

\draw (10,0) node {$\simeq$};

\draw (11,0) node {\scalebox{3}{$\int$}};
\draw[
        decoration={markings, mark=at position 0 with {\arrow{>}}},
        decoration={markings, mark=at position 0.5 with {\arrow{>}}},
        postaction={decorate}
        ] (12,0) circle (0.5cm);
\draw (11.5,0) node[anchor=east] {$\cC$};
\draw (12.5,0) node[anchor=west] {$\cD$};
\draw ({12+0.5*cos(60)},{0.5*sin(60)}) node {$\bullet$};
\draw ({12+0.5*cos(-60)},{0.5*sin(-60)}) node {$\bullet$};
\draw ({12+0.5*cos(60)+0.2},{0.5*sin(60)}) node[anchor=south] {$\cP$};
\draw ({12+0.5*cos(-60)+0.2},{0.5*sin(-60)}) node[anchor=north] {$\cQ$};


\draw (-1.5,-2) node {$\simeq$};

\draw (-0.5,-2) node {\scalebox{3}{$\int$}};
\draw[
        decoration={markings, mark=at position 0 with {\arrow{>}}},
        decoration={markings, mark=at position 0.5 with {\arrow{>}}},
        postaction={decorate}
        ] (0.5,-2) circle (0.5cm);
\draw (0,-2) node[anchor=east] {$\cC$};
\draw (1,-2) node[anchor=west] {$\cD$};
\draw ({0.5+0.5*cos(120)},{-2+0.5*sin(120)}) node {$\bullet$};
\draw ({0.5+0.5*cos(-120)},{-2+0.5*sin(-120)}) node {$\bullet$};
\draw ({0.5+0.5*cos(120)-0.2},{-2+0.5*sin(120)}) node[anchor=south] {$\cP$};
\draw ({0.5+0.5*cos(-120)-0.2},{-2+0.5*sin(-120)}) node[anchor=north] {$\cQ$};

\draw (2,-2) node {$\simeq$};

\draw (3,-2) node {\scalebox{3}{$\int$}};
\draw[
        decoration={markings, mark=at position 0 with {\arrow{>}}},
        postaction={decorate}
        ] (5,-2) circle (0.5cm);
\draw ({4.5+0.025},-2) node[anchor=east] {$(\cP \otimes_\cC \! \cQ)$};
\draw (4.5,-2) node {$\bullet$};
\draw (5.5,-2) node[anchor=west] {$\cD$};

\draw (6.5,-2) node {$\simeq$};

\draw (7.5,-2) node {\scalebox{3}{$\int$}};
\draw[
        decoration={markings, mark=at position 0 with {\arrow{>}}},
        postaction={decorate}
        ] (8.5,-2) circle (0.5cm);
\draw (8,-2) node {$\bullet$};
\draw (8,-2) node[anchor=east] {$\cD$};
\draw (9,-2) node[anchor=west] {$\cD$};

\draw (10,-2) node {$\simeq$};

\draw (11,-2) node {\scalebox{3}{$\int$}};
\draw[
        decoration={markings, mark=at position 0 with {\arrow{>}}},
        postaction={decorate}
        ] (12,-2) circle (0.5cm);
\draw (12.5,-2) node[anchor=west] {$\cD$};

\end{tikzpicture}
\]
in $\cV$ among the indicated values of ``$\cB$-framed'' factorization homology \cite{AFT-fact}, a generalization of factorization homology that allows distinct strata to be labeled by distinct algebraic data.\footnote{It would require a nontrivial amount of rigorously merge the formalism of $\cB$-framed factorization homology developed in \cite{AFT-fact} with the formalism of enriched factorization homology developed in \cite{AMR-fact}.  On the other hand, note that this assertion is proved in \cite{BluMan-loc} (in the ``set of objects'' and spectrally-enriched setting).}
\end{remark}

\begin{remark}
\label{functor.gives.bimod}
An arbitrary functor
\[
\cC
\longra
\cD
\]
of flagged $\cV$-enriched $\infty$-categories induces a $(\cC,\cD)$-bimodule
\[
\cC^\op \times \cD
\longra
\cD^\op \times \cD
\xra{\ulhom_\cD(-,-)}
\cV
~.
\]
This allows us to consider the $\infty$-category $\fCat(\cV)$ of flagged $\cV$-enriched $\infty$-categories as a full subcategory of a certain Morita $\infty$-category.\footnote{It is indeed merely a condition for a bimodule to be represented by an actual functor: this is the condition that its adjunct factors through the Yoneda embedding, a monomorphism.}  Thus, bimodules can be considered as a generalization of ordinary functors, and in particular it is meaningful to assert that a morphism in $\fCat(\cV)$ is a Morita equivalence.

Let us note two examples of this phenomenon of particular interest.
\begin{enumerate}

\item Let $A \in \Alg(\Spectra)$ be a ring spectrum, and write $\fB A$ for its corresponding one-object flagged $\cV$-enriched $\infty$-category.  Then there is a canonical spectrally-enriched functor
\[
\fB A
\longra
\Perf_A
~,
\]
which takes the unique object of $\fB A$ to the object $A \in \Perf_A$ and induces an equivalence $A \xra{\sim} \ulhom_{\Perf_A}(A,A)$ on hom-objects.  This is a Morita equivalence because $\Perf_A$ is the \textit{stable idempotent-completion} of $\fB A$, and so the restriction
\[
\ul{\Fun}(\fB A,\Spectra)
\xlongla{\sim}
\ul{\Fun}(\Perf_A,\Spectra)
\simeq \Fun^{\sf ex}(\Perf_A,\Spectra)
\]
is an equivalence (where $\ul{\Fun}$ denotes the $\infty$-category of spectrally-enriched functors and $\Fun^{\sf ex}$ denotes the $\infty$-category of exact functors).

\item \label{compln.is.mor.invt} Let $\cC \in \fCat(\cV)$ be any flagged $\cV$-enriched $\infty$-category, and write $\ol{\cC} \in \Cat(\cV)$ for its completion (see \Cref{complete.V.enr.cats}).  Then, the canonical functor
\[
\cC
\longra
\ol{\cC}
\]
is a Morita equivalence: the restriction
\[
\Mod_\cC
:=
\ul{\Fun}(\cC,\cV)
\xlongla{\sim}
\ul{\Fun}(\ol{\cC},\cV)
=:
\Mod_{\ol{\cC}}
\]
is an equivalence since $\cV$ itself is complete.\footnote{By passing from $\cV$ to the $\infty$-category of presheaves on $\cV$ equipped with Day convolution (and noting that the Yoneda embedding is a monoidal monomorphism), one can assume without loss of generality that $\cV$ is presentably (and hence closed) monoidal.}
\end{enumerate}
In particular, combining example \Cref{compln.is.mor.invt} with \Cref{morita.invce} implies that our definition of $\THH$ (for spectrally-enriched $\infty$-categories) agrees with those defined in terms of a set of objects, as e.g.\! in \cite{BluMan-loc}.
\end{remark}


\begin{remark}
\label{fgt.cyclo.is.conservative}
The forgetful functor
\[
\Cyclo(\Spectra)
\xlongra{\fgt}
\Spectra
\]
is conservative, as both functors in the composite
\[
\Cyclo(\Spectra)
\xlongra{\fgt}
\Fun(\BT,\Spectra)
\xlongra{\fgt}
\Spectra
\]
are easily seen to be conservative.  In particular, for any ring spectrum $A \in \Alg(\Spectra)$, the canonical functor
\[
\fB A
\longra
\Perf_A
\]
of \Cref{functor.gives.bimod} induces an equivalence
\[ 
\THH(\fB A)
\xlongra{\sim}
\THH(\Perf_A)
\]
of cyclotomic spectra, and thereafter an equivalence
\[
\TC(\fB A)
\xlongra{\sim}
\TC(\Perf_A)
~.
\]
More generally, for any $\cC \in \fCat(\Spectra)$, the Yoneda embedding $\cC \ra \Perf_\cC$ (into the stable idempotent-completion) is a Morita equivalence.  Moreover, Morita equivalences between stable idempotent-complete $\infty$-categories always arise from actual functors.  Thus, \textit{all} Morita equivalences in $\fCat(\Spectra)$ can be realized as zigzags of enriched functors; hence, $\THH$ takes them to equivalences of cyclotomic spectra, and $\TC$ takes them to equivalences of spectra.
\end{remark}

\subsection{Outline}

This paper is entirely self-contained: over the course of our proofs, we recall everything that we need from the supporting papers \cite{AMR-cyclo,AMR-fact}.  It is organized as follows.
\begin{itemize}

\item In \Cref{part.cyclo.str}, we prove \Cref{THH.is.cyclo} (the cyclotomic structure on $\THH$).  We begin with a review of stratified 1-manifolds in \Cref{section.roadmap.stratified.1.mflds}, and of enriched factorization homology in \Cref{section.roadmap.fact.hlgy}.  We then study the unstable cyclotomic structure on spatially-enriched factorization homology over the circle in \Cref{section.roadmap.unst.cyclo.str}, in a way that anticipates its stable analog.  Finally, in \Cref{section.roadmap.cyclo.str} we construct the cyclotomic structure on spectrally-enriched factorization homology over the circle (i.e.\! $\THH$).

\item In \Cref{part.cyclo.trace}, we deduce \Cref{cyclo.trace} (the cyclotomic trace).  In \Cref{section.cyclo.trace}, we review the unstable cyclotomic trace.  Then, in \Cref{section.unstable.cyclo.trace} we use this to obtain the cyclotomic trace.

\end{itemize}

\subsection{Notation and conventions}
\label{subsection.notation.and.conventions}

\begin{enumerate}

\item \catconventions \inftytwoconventions

\item \functorconventions

\item \circconventions


\item \spacescatsspectraconventions


\item \fibrationconventions

\item\label{efibconventions} \efibconventions

\item \cavalieraboutbicompleteness

\item \codenamesfortrace

\end{enumerate}

\subsection{Acknowledgments}

We thank Andrew Blumberg for numerous helpful conversations related to this project.  We also thank Hood Chatham and Eric Peterson for their TeX support.

\acksupport

\stopcontents[sections]

\part{The cyclotomic structure on $\THH$}
\label{part.cyclo.str}

In this part of the paper, we prove \Cref{THH.is.cyclo} modulo certain technical details which are proved in the papers \cite{AMR-cyclo,AMR-fact}, giving precise and specific references for all such details as we proceed.

To begin, in \Cref{section.roadmap.stratified.1.mflds} we describe the essential features of \textit{stratified 1-manifolds}, which parametrize our entire story.

Although $\THH$ is given specifically by \textit{spectrally}-enriched factorization homology over the circle, we will find its cyclotomic structure as a consequence of a more primitive structure on \textit{spaces}-enriched factorization homology over the circle.  Thus, in \Cref{section.roadmap.fact.hlgy} we discuss $\cV$-enriched factorization homology over the circle for an arbitrary symmetric monoidal $\infty$-category $(\cV,\boxtimes)$, which we denote by
\[ \THHV~. \]
We write
\[ \THHVx \]
when we specialize to the case that $(\cV,\boxtimes) = (\cV,\times)$ is \textit{cartesian} symmetric monoidal, and we simply write
\[ \THH := \THH_\Spectra \]
when we specialize to the particular case that $(\cV,\boxtimes ) = (\Spectra,\otimes)$ is the $\infty$-category of spectra (equipped with its smash product symmetric monoidal structure).

In \Cref{section.roadmap.unst.cyclo.str} \and \Cref{section.roadmap.cyclo.str}, we proceed to incorporate additional symmetries via those of 1-manifolds.  There exists a canonical $\TT$-action on $\THHV$ for any symmetric monoidal enriching $\infty$-category $\cV$, and $\THHVx$ carries additional contravariant functoriality for covering maps of circles via diagonal maps in $\cV$; we distill this state of affairs in \Cref{section.roadmap.unst.cyclo.str}, taking an expositional detour to indicate how the unstable cyclotomic structure arises in order to motivate the analogous construction in the stable case.  Finally, in \Cref{section.roadmap.cyclo.str} we show how the diagonal maps for spaces induce Tate diagonal maps for spectra, which in turn endow $\THH$ with its cyclotomic structure.

\section{Stratified 1-manifolds}
\label{section.roadmap.stratified.1.mflds}

Factorization homology in dimension 1 is parametrized by a certain $\infty$-category
\[ \M \]
of \bit{stratified 1-manifolds}.  We describe its essential features in \Cref{subsection.strat.mflds}, and in \Cref{subsection.disk.strats} we show how it gives rise to various indexing categories of classical interest.  The details omitted from this section are contained in [\Cref{fact:section.strat.mflds}].

\subsection{Stratified 1-manifolds}
\label{subsection.strat.mflds}

The objects of $\M$ are simply finite disjoint unions of framed circles and finite directed connected graphs and framed circles. 
For any stratified 1-manifold $M \in \M$, we write $M^{(0)}$ and $M^{(1)}$ for its 0- and 1-dimensional strata, which we'll often identify with their underlying spaces. 
The full subcategory
\[ \D \subset \M \]
of \bit{disk-stratified} 1-manifolds then consists of those objects whose strata are both disjoint unions of disks; equivalently, an object of $\M$ is disk-stratified just when it doesn't contain any (unstratified) circles among its connected components.

While the morphisms in $\M$ are somewhat more complicated to describe, it will suffice to highlight two distinguished classes thereof.
\begin{itemize}
\item The \bit{refinement morphisms} correspond \textit{covariantly} to homeomorphisms on underlying topological spaces, but where in passing from the source to the target we're allowed to forget some of the marked points (while the framings on the 1-dimensional strata must coincide).  
\item The \bit{creation morphisms} correspond \textit{contravariantly} to surjective proper constructible bundles. 
\end{itemize}

We will generally denote an arbitrary object of $\M$ by $M$; we will write $S$ to specifically denote a (smooth) circle, and $R$ to denote an object of $\D$ (which will frequently be the source of a refinement morphism).

\subsection{Disk-refinements}
\label{subsection.disk.strats}

Factorization homology over a stratified 1-manifold $M \in \M$ will be given by a colimit over the $\infty$-category $\D(M)$ of \bit{disk-refinements} of $M$, an object of which is simply a disk-stratified 1-manifold $R$ equipped with a refinement morphism $R \ra M$: writing
\[
\totalD
\subset 
\Ar(\M)
\]
for the full subcategory of the arrow $\infty$-category of $\M$ on the refinement morphisms with disk-stratified source, this is given by the fiber
\[
\begin{tikzcd}
\D(M)
\arrow{r}
\arrow{d}
&
\totalD
\arrow{d}{t}
\\
\{ M \}
\arrow{r}
&
\M
\end{tikzcd}
~.
\]
In the special case that we take $M = S^1$ to be the circle, we have an identification
\[
\D(S^1)
\simeq
\para
\]
with the \bit{paracyclic indexing category}.  This carries a left action of the circle group $\TT$ coming from the functoriality of the construction $\D(-)$ for equivalences in $\M$, and its quotient thereby admits the further identification
\[
(\para)_{\htpy \TT}
\simeq
\cyclic
\]
with the \bit{cyclic indexing category}.  In fact, this -- as well as a larger category of interest -- can be recovered directly as follows.  Let us respectively write
\[
\BT
\longrsurjmono
\BW
\xlonghookra{\ff}
\M
\]
for the full subgroupoid and the full subcategory on the object $S^1 \in \M$; the monoid $\WW$ can be identified as the semidirect product
\[
\WW
\simeq
\TT
\rtimes
\Nx
\]
classified by the right action of the monoid $\Nx$ on the group $\TT$ by positive-degree endomorphisms.\footnote{The notation $\WW$ stems from the fact that this will keep track of Frobenius and Verschiebung operators, along the lines of \cite{HessMad-Witt}.}  Then, we have a diagram
\begin{equation}
\label{pullbacks.from.D.ref.over.M}
\begin{tikzcd}
\para
\arrow[two heads]{r}
\arrow{d}
&
\cyclic
\arrow[hook, two heads]{r}
\arrow{d}
&
\epicyc
\arrow[hook]{r}{\ff}
\arrow{d}
&
\totalD
\arrow{d}{t}
\\
\{S^1\}
\arrow[two heads]{r}
&
\BT
\arrow[hook, two heads]{r}
\arrow{d}
&
\BW
\arrow[hook]{r}[swap]{\ff}
\arrow{d}
&
\M
\\
&
\pt
\arrow[hook, two heads]{r}
&
\BN
\end{tikzcd}
~,
\end{equation}
in which
\begin{itemize}
\item all squares are pullbacks,
\item among the upper vertical maps, all but the rightmost are cocartesian fibrations,\footnote{The rightmost map $\totalD \ra \M$ is not a cocartesian fibration, but it is when restricted to a fairly large subcategory [\Cref{fact:functoriality.for.D.of}].}
\item both lower vertical maps are right fibrations, and
\item $\epicyc$ denotes the \bit{epicyclic indexing category}.
\end{itemize}
In particular, the epicyclic category is given by the left-lax quotient
\[
\epicyc
\simeq
(\para)_{\htpy^\llax \WW}
 ~ ,
\]
where the left action of the monoid $\WW$ arises from the contravariant functoriality of the construction $\D(-)$ for the positive-degree self-coverings of the circle: disk-refinements can be pulled back along covering maps.

Whereas the paracyclic category $\para$ might be thought of as being given by disk-refinements of \textit{the} circle, the cyclic and epicyclic categories $\cyclic$ and $\epicyc$ might be thought of as being given by disk-refinements of \textit{a} circle: the point being that in the morphisms in these latter two categories the underlying circle need not remain fixed, but can change either by an automorphism or an endomorphism in $\M$, respectively.

\section{Factorization homology}
\label{section.roadmap.fact.hlgy}

In this section we define
\[
\THHV(\cC)
~,
\]
the \bit{enriched factorization homology over the circle} of an arbitrary (flagged) $\cV$-enriched $\infty$-category $\cC$.  This definition is given in \Cref{subsection.enr.fact.hlgy}.  It relies on a \textit{categorified} form of factorization homology, which we describe in \Cref{subsection.catfied.fact.hlgy}.  We begin in \Cref{subsection.enr.infty.cats} by recalling the formalism of $\cV$-enriched $\infty$-categories laid out in \cite{GH-enr}.  The details omitted from this section are contained in [\Cref{fact:section.categorified.fact.hlgy} \and \Cref{fact:section.enriched.fact.hlgy}].

\subsection{Enriched $\infty$-categories}
\label{subsection.enr.infty.cats}

Let $\cV := (\cV,\boxtimes)$ be a monoidal $\infty$-category, and write
\[ \bDelta^\op \xra{\BV} \Cat \]
for its \textit{monoidal deloop} -- that is, its bar construction.  This is a \textit{category object} in $\Cat$, i.e.\! a simplicial object satisfying the Segal condition.  Then, a \bit{$\cV$-enriched $\infty$-category} $\cC$ can be specified as follows.  First, we must specify its \textit{underlying $\infty$-groupoid}, an object $\iC \in \Spaces$.  Let us write
\[ \bDelta^\op \xra{\cdiC} \Spaces \]
for the ``codiscrete category object'' on $\iC$: the category object in $\Spaces \subset \Cat$ whose space of objects is $\iC$ and whose hom-spaces are all contractible.  Then, $\cC$ is specified by its \bit{enriched hom functor}, a right-lax functor
\begin{equation}
\label{enr.hom.functor.for.C}
\cdiC
\xra{\ulhom_\cC}
\BV
\end{equation}
of category objects in $\Cat$.  In simplicial degree $n$, this is given by
\[ \begin{tikzcd}[row sep=0cm]
(\iC)^{\times (n+1)}
\arrow{r}
&
\cV^{\times n}
\\
\rotatebox{90}{$\in$}
&
\rotatebox{90}{$\in$}
\\
(C_0,\ldots,C_n)
\arrow[mapsto]{r}
&
( \ulhom_\cC(C_0,C_1) , \ldots , \ulhom_\cC(C_{n-1},C_n) )
\end{tikzcd}
~,
\]
and the right-laxness determines the categorical structure maps; for instance, restricting to the indicated morphisms in $\bDelta^\op$, this specifies the diagram
\[ \begin{tikzcd}
{[2]^\circ}
\arrow{d}[swap]{\delta_1}
&
(\iC)^{\times 3}
\arrow{r}
\arrow{d}
&
\cV^{\times 2}
\arrow{d}
\\
{[1]^\circ}
&
(\iC)^{\times 2}
\arrow{r}[pos=0.3, transform canvas={yshift=0.4cm}]{\rotatebox{-135}{$\Rightarrow$}}[swap, pos=0.3, transform canvas={yshift=-0.4cm}]{\rotatebox{135}{$\Rightarrow$}}
&
\cV
\\
{[0]^\circ}
\arrow{u}{\sigma_0}
&
\iC
\arrow{r}
\arrow{u}
&
\pt
\arrow{u}
\\
\bDelta^\op
&
\cdiC
\arrow{r}{\ulhom_\cC}
&
\BV
\end{tikzcd}
~,
\]
in which the upper square selects the composition maps
\begin{equation}
\label{composition.in.enr.cat}
\ulhom_\cC(C_0,C_1) \boxtimes \ulhom_\cC(C_1,C_2)
\longra
\ulhom_\cC(C_0,C_2)
\end{equation}
while the lower square selects the unit maps
\begin{equation}
\label{unit.in.enr.cat}
\uno_\cV
\longra
\ulhom_\cC(C_0,C_0)
~.
\end{equation}

In fact, the above data merely specify a \bit{flagged} $\cV$-enriched $\infty$-category; these assemble into an $\infty$-category
\[
\fCat(\cV)
~.
\]
Among these, one can demand a \textit{completeness} condition: there is an internally-defined ``maximal subgroupoid'' which receives a $\pi_0$-surjective map from the underlying $\infty$-groupoid $\iC$, and the condition is that this map be an equivalence.\footnote{This is analogous to the completeness condition for Segal spaces: the map from a Segal space to its completion is always a $\pi_0$-surjection in simplicial degree 0, and the Segal space is equivalent data to its completion (or equivalently its corresponding $\infty$-category) along with this $\pi_0$-surjection.}  These define a full subcategory
\[
\Cat(\cV)
\subset
\fCat(\cV)
~,
\]
which admits a left adjoint when $\cV$ is presentably monoidal.  However, as mentioned in \Cref{flagged}, we will have no reason to restrict our attention to this subcategory, and indeed when we take $\iC = \pt$ we recover the definition of the monoidal deloop $\cC = \fB A$ of an algebra object $A \in \Alg(\cV)$ (and thereafter the classical definition of $\THH(A)$).

\subsection{Categorified factorization homology}
\label{subsection.catfied.fact.hlgy}

Now, we can take factorization homology \textit{of category objects} in $\Cat$.  In fact, we will only need this for disk-stratified 1-manifolds, and this is particularly easy to describe: for a disk-stratified 1-manifold $R \in \D$ and a category object $\cY \in \Fun(\bDelta^\op,\Cat)$, the factorization homology
\[ \int_R \cY \in \Cat \]
is given by a limit of copies of $\cY_{|[0]^\circ}$ and $\cY_{|[1]^\circ}$: one copy of $\cY_{|[0]^\circ}$ for each vertex, one copy of $\cY_{|[1]^\circ}$ for each edge, and with structure maps $s = \delta_1$ and $t = \delta_0$ determined by the incidence data of the directed graph $R$.\footnote{The composition in $\cY$ plays a role in the functoriality of this construction, and in particular in the definition of the factorization homology $\int_M \cY$ when $M \in \M$ is no longer disk-stratified.}  This is suitably functorial for right-lax functors of category objects, so that a $\cV$-enriched $\infty$-category $\cC$ -- that is, its enriched hom functor \Cref{enr.hom.functor.for.C} -- determines a morphism
\begin{equation}
\label{categorified.fact.hlgy.over.D}
\begin{tikzcd}[column sep=2cm]
\int_{|\D} \cdiC
\arrow{r}{\int_{|\D} \ulhom_\cC}
\arrow{d}
&
\int_{|\D} \BV
\arrow{ld}
\\
\D
\end{tikzcd}
\end{equation}
in $\Cat_{\cocart/\D}$.\footnote{The morphism \Cref{categorified.fact.hlgy.over.D} does not lie in $\coCart_{\D}$, but it does preserve cocartesian morphisms over a fairly substantial class of morphisms in $\D$ [\Cref{fact:obs.cocart.over.cls.and.some.cr}].}  In particular, pulling back along the composite
\[
\para
\longra
\totalD
\xlongra{s}
\D
~,
\]
we obtain a likewise morphism
\begin{equation}
\label{categorified.fact.hlgy.over.para}
\begin{tikzcd}[column sep=2cm]
\int_{|\para} \cdiC
\arrow{r}{\int_{|\para} \ulhom_\cC}
\arrow{d}
&
\int_{|\para} \BV
\arrow{ld}
\\
\para
\end{tikzcd}
\end{equation}
in $\Cat_{\cocart/\para}$.  Over a disk-refinement $(R \ra S^1) \in \para$, an object of its source is simply given by a labeling of its set of vertices $R^{(0)}$ by objects of the $\infty$-groupoid
\[ \iC =: \cdiC_{|[0]^\circ} \]
(since the morphism-data in $\cdiC$ are canonically determined), while an object of its target is simply given by a labeling of its set of edges $R^{(1)}$ by objects of the $\infty$-category
\[ \cV =: \BV_{|[1]^\circ} \]
(since the object-data in $\BV$ are canonically determined).  In other words, over this object of $\para$ the map \Cref{categorified.fact.hlgy.over.para} restricts to a map
\[
(\iC)^{\times R^{(0)}}
\longra
\cV^{\times R^{(1)}}
~.
\]
Of course, this is given by nothing other than the enriched hom functor of $\cC$. 

\subsection{Enriched factorization homology}
\label{subsection.enr.fact.hlgy}

When $\cV$ is additionally \textit{symmetric} monoidal, there exists a ``tensor everything together'' functor
\[
\int_{|\para} \BV
\xlongra{\bigboxtimes}
\cV
~.\footnote{In fact, in order to guarantee the existence of this functor, it suffices (definitionally) for $\cV$ to be merely cyclically monoidal.}
\]
Then, we define the \bit{enriched factorization homology} of $\cC$ over $S^1$ 
by the formula
\begin{equation}
\label{define.THH}
\THHV(\cC)
:=
\int_{S^1} \cC
:=
\colim
\left(
\int_{|\para} \cdiC
\xra{\int_{|\para} \ulhom_\cC}
\int_{|\para} \BV
\xlongra{\bigboxtimes}
\cV
\right)
~.
\end{equation}

Note that this indeed precisely rigorizes the heuristic definition of \Cref{subsection.cocycles}.  First of all, this colimit is indexed over the data of collections of points $p_1,\ldots,p_n \in S^1$ labeled by objects $C_1,\ldots,C_n \in \cC$, with value the monoidal product
\[
\underset{i \in \ZZ/n}{\standaloneboxtimes} \ulhom_\cC(C_i,C_{i+1})
~.
\]
Moreover, its structure maps are given either by disappearances of points (i.e.\! refinement morphisms among disk-refinements of $S^1$), which are taken to composition maps (as map \Cref{composition.in.enr.cat}), or by anticollisions of points (i.e.\! refinement morphisms among disk-refinements of $S^1$), which are taken to unit maps (as map \Cref{unit.in.enr.cat}).

This definition of enriched factorization homology has the key feature of being optimized for the isolation of the various specific properties of the enriching $\infty$-category $\cV$: all natural operations thereon have only to do with the factorization homology of $\BV$, while the factorization homology of $\cdiC$ simply comes along for the ride.

\section{The unstable cyclotomic structure on $\THHVx$}
\label{section.roadmap.unst.cyclo.str}

In this section, we obtain the \bit{unstable cyclotomic structure}
\[
\Cyclo^\htpy
\lacts
\THHVx
\]
in the case that $(\cV,\boxtimes) = (\cV,\times)$ is \textit{cartesian} symmetric monoidal.  In fact, this structure is nothing but an action of the monoid $\WW$, which we will obtain in \Cref{subsection.cartesian.fact.hlgy.with.W.action} after obtaining its underlying $\TT$-action in \Cref{subsection.fact.hlgy.with.T.action}.  These are both special cases of more general constructions: \Cref{subsection.fact.hlgy.with.T.action} is based on material in [\Cref{fact:section.enriched.fact.hlgy}] while \Cref{subsection.cartesian.fact.hlgy.with.W.action} is based on material in [\Cref{fact:section.cartesian.enr.fact.hlgy}], and all details omitted from these two subsections can be found there.

In order to motivate our construction of the cyclotomic structure on $\THH := \THH_\Spectra$, we will rederive the unstable cyclotomic structure on $\THHVx$ as a consequence of the \bit{diagonal package}, which packages the diagonal maps that exist in $\cV$.  We accomplish this as follows.  First, in \Cref{subsection.diagonal.package}, we describe the diagonal package.  Then, in \Cref{subsection.unstable.cyclo.objects}, we reformulate our definition of unstable cyclotomic objects in similar terms.  Finally, in \Cref{subsection.unstable.cyclo.str}, we describe the way in which the diagonal package induces the unstable cyclotomic structure on $\THHVx$.

The diagonal package will be a global section of a certain cartesian fibration over $\BN$.  So before describing the diagonal package, we study this cartesian fibration: we provide a heuristic description thereof in \Cref{subsection.local.str.of.Prod}, and we give a careful analysis (as well as a useful generalization) thereof in \Cref{subsection.cart.mdrmy.via.cocart.mdrmy}.

Given the underlying $\TT$-action on $\THHVx$, the extension to a $\WW$-action -- that is, the data of its unstable cyclotomic structure -- is manifestly equivalent data to a family of \textit{unstable cyclotomic structure maps}
\begin{equation}
\label{unst.cyclo.str.map}
\THHVx
\longra
\left( \THHVx \right)^{\htpy \Cyclic_r}
\end{equation}
for all $r \in \Nx$, which are compatible in the sense that they compose in the evident way.\footnote{Note the similarity with our description of a cyclotomic structure on a $\TT$-spectrum in \Cref{tau.action.only.left.lax}.}  For expository purposes, in \Cref{subsection.unwind.unst.cyclo.str.map} we show explicitly how the unstable cyclotomic structure maps \Cref{unst.cyclo.str.map} arise from diagonal maps in $\cV$.

\subsection{The $\TT$-action on $\THHV$}
\label{subsection.fact.hlgy.with.T.action}

For any symmetric monoidal $\infty$-category $\cV$, we can obtain a left $\TT$-action on $\THHV$ by a direct extension of formula \Cref{define.THH}.  Namely, rather than pulling the map \Cref{categorified.fact.hlgy.over.D} all the way back along the composite
\[
\para
\longra
\cyclic
\longra
\epicyc
\longra
\totalD
\xlongra{s}
\D
\]
(recall diagram \Cref{pullbacks.from.D.ref.over.M}), we can instead pull it back only to the cyclic category $\cyclic$.  We again have a ``tensor everything together'' functor
\begin{equation}
\label{tensor.everything.together.over.cyclic}
\int_{|\cyclic} \BV
\xlongra{\bigboxtimes}
\cV
~.
\end{equation}
So, we obtain the desired action by taking the left Kan extension
\begin{equation}
\label{fact.hlgy.with.T.action}
\begin{tikzcd}[column sep=2cm, row sep=1.5cm]
\int_{|\cyclic} \cdiC
\arrow{r}{\int_{|\cyclic} \ulhom_\cC}
\arrow{d}
&
\int_{|\cyclic} \BV
\arrow{ld}
\arrow{r}{\boxtimesfortikz}
&
\cV
\\
\cyclic
\arrow{d}
\\
\BT
\arrow[dashed, bend right=10]{rruu}[swap]{\TT \lacts \THHV(\cC)}
\end{tikzcd}
~;
\end{equation}
as the two vertical functors are cocartesian fibrations then so is their composite, and hence this simply reduces to a fiberwise colimit.\footnote{As the fiber of $\cyclic$ over each point in $\BT$ is a copy of $\para$, over each point in $\BT$ we're taking a colimit over the fiber
\[ \int_{|\para} \cdiC ~ ; \]
thus, the underlying object is indeed $\THHV(\cC)$, as defined in formula \Cref{define.THH}.}

\subsection{The $\WW$-action on $\THHVx$}
\label{subsection.cartesian.fact.hlgy.with.W.action}

Now, one might hope to extend the construction of \Cref{subsection.fact.hlgy.with.T.action} yet further, from the cyclic category $\cyclic$ to the epicyclic category $\epicyc$.  However, one immediately runs into a problem: there is only a ``tensor everything together'' functor
\[
\int_{|\epicyc} \BV
\xlongra{\bigboxtimes}
\cV
\]
when $\cV$ is \textit{cartesian} symmetric monoidal;\footnote{In fact, it suffices (definitionally) for $\cV$ to be an augmented cyclically monoidal and cyclically comonoidal bialgebra object in $\Cat$: this is precisely the data of a suitable functor generalizing the one labeled $\prod$ in diagram \Cref{fact.hlgy.with.W.action}.} a simple case in which this issue arises is e.g.\! for the morphism
\begin{equation}
\label{pic.of.morphism.in.int.BV.over.three.fold.creation}
\begin{tikzpicture}

\draw[
        decoration={markings, mark=at position 0.5 with {\arrow{>}}},
        postaction={decorate}
        ] (-0.5,0) circle (0.5cm);

\draw (0,0) node {$\bullet$};
\draw (-1,0) node[anchor=east] {$V$};

\draw (1,0) node {$\longra$};

\draw[
        decoration={markings, mark=at position 1/6 with {\arrow{>}}},
        decoration={markings, mark=at position 0.5 with {\arrow{>}}},
        decoration={markings, mark=at position 5/6 with {\arrow{>}}},
        postaction={decorate}
        ] (4,0) circle (1.5cm);

\draw (5.5,0) node {$\bullet$};
\draw ({4+1.5*cos(120)},{0+1.5*sin(120)}) node {$\bullet$};
\draw ({4+1.5*cos(240)},{0+1.5*sin(240)}) node {$\bullet$};

\draw ({2.5-0.05},0) node[anchor=east] {$V$};
\draw ({4+1.5*cos(60)+0.25},{0+1.5*sin(60)+0.05}) node[anchor=south] {$V$};
\draw ({4+1.5*cos(-60)+0.25},{0+1.5*sin(-60)-0.05}) node[anchor=north] {$V$};

\end{tikzpicture}
\end{equation}
in $\int_{|\epicyc} \BV$, which must be sent to a natural morphism
\[
V
\longra
V^{\boxtimes 3}
\]
in $\cV$.  However, in this case, we can indeed obtain a $\WW$-action on $\THHVx(\cC)$ by taking the left Kan extension
\begin{equation}
\label{fact.hlgy.with.W.action}
\begin{tikzcd}[column sep=2cm, row sep=1.5cm]
\int_{|\epicyc} \cdiC
\arrow{r}{\int_{|\epicyc} \ulhom_\cC}
\arrow{d}
&
\int_{|\epicyc} \BV
\arrow{ld}
\arrow{r}{\prod}
&
\cV
\\
\epicyc
\arrow{d}
\\
\BW
\arrow[dashed, bend right=10]{rruu}[swap]{\WW \lacts \THHVx(\cC)}
\end{tikzcd}
~,
\end{equation}
again a fiberwise colimit.\footnote{An explicit description of the functoriality of this fiberwise colimit is given in \Cref{subsection.unwind.unst.cyclo.str.map}.}  This will turn out to be precisely the data of the \textit{unstable cyclotomic structure} on $\THHVx(\cC)$.

\subsection{The local structure of the functor $\prod$}
\label{subsection.local.str.of.Prod}

To see more clearly what is going on, it is helpful to work ``locally'' over $\BN$.  First of all, to specify the functor $\prod$ in diagram \Cref{fact.hlgy.with.W.action} is equivalently to specify a functor
\[ \begin{tikzcd}
&
\int_{|\epicyc} \BV
\arrow[dashed]{r}{\prod}
\arrow{ld}
&
\ul{\cV}
\arrow{dd}
&
\hspace{-1.3cm} := \cV \times \BN
\\
\epicyc
\arrow{d}
\\
\BW
\arrow{rr}
&
&
\BN
\end{tikzcd}
\]
over $\BN$.  This, in turn, we can consider as a global section
\begin{equation}
\label{global.section.of.relfun}
\begin{tikzcd}[row sep=1.5cm]
\Fun^\rel_{/\BN} \left( \int_{|\epicyc} \BV , \ul{\cV} \right)
\arrow{d}
\\
\BN
\arrow[dashed, bend left]{u}{\prod}
\end{tikzcd}
\end{equation}
of the \textit{relative functor $\infty$-category} (recall item \Cref{efibconventions} of \Cref{subsection.notation.and.conventions}).  Now, the fiber of the downwards functor in diagram \Cref{global.section.of.relfun}
over the unique object of $\BN$ is given by
\begin{equation}
\label{fiber.of.cartesian.fibn.over.BN}
\Fun \left( \int_{|\cyclic} \BV , \cV \right) ~.
\end{equation}
Moreover, as we will see in \Cref{subsection.cart.mdrmy.via.cocart.mdrmy}, a core computation reveals that this downwards functor is in fact a cartesian fibration: over a morphism $[1] \xra{r} \BN$, the cartesian monodromy is given heuristically by the functor
\begin{equation}
\label{cartesian.mdrmy.by.h}
\begin{tikzcd}[row sep=0cm]
\Fun \left( \int_{|\cyclic} \BV , \cV \right)
&
\Fun \left( \int_{|\cyclic} \BV , \cV \right)
\arrow{l}
\\
\rotatebox{90}{$\in$}
&
\rotatebox{90}{$\in$}
\\
\left(
\scalebox{0.5}{
\begin{tikzpicture}

\pgfmathsetmacro{\R}{1.5}

\draw[
        decoration={markings, mark=at position 15/360 with {\arrow{>}}},
        decoration={markings, mark=at position 45/360 with {\arrow{>}}},
        decoration={markings, mark=at position 75/360 with {\arrow{>}}},
        decoration={markings, mark=at position -15/360 with {\arrow{>}}},
        decoration={markings, mark=at position -45/360 with {\arrow{>}}},
        postaction={decorate}
        ] (0,0) circle (\R cm);

\draw ({\R*cos(0)},{\R*sin(0)}) node {$\bullet$};
\draw ({\R*cos(30)},{\R*sin(30)}) node {$\bullet$};
\draw ({\R*cos(60)},{\R*sin(60)}) node {$\bullet$};
\draw ({\R*cos(-30)},{\R*sin(-30)}) node {$\bullet$};

\draw ({\R*(cos(15)+0.05)},{\R*(sin(15)-0.05)}) node[anchor=west] {$V_1$};
\draw ({\R*cos(45)},{\R*(sin(45)+0.1)}) node[anchor=west] {$V_2$};
\draw ({\R*(cos(75)+0.1)},{\R*sin(75)}) node[anchor=south] {\rotatebox{-15}{$\cdots$}};

\draw ({\R*(cos(-15)+0.05)},{\R*(sin(-15)-0.1)}) node[anchor=west] {$V_n$};
\draw ({\R*(cos(-45)-0.1)},{\R*(sin(-45)-0.15)}) node[anchor=west] {\rotatebox{45}{$\cdots$}};

\end{tikzpicture}
}
\longmapsto
F \left(
\scalebox{0.5}{
\begin{tikzpicture}

\pgfmathsetmacro{\R}{3}

\draw[
        decoration={markings, mark=at position 15/720 with {\arrow{>}}},
        decoration={markings, mark=at position 45/720 with {\arrow{>}}},
        decoration={markings, mark=at position 75/720 with {\arrow{>}}},
        decoration={markings, mark=at position -15/720 with {\arrow{>}}},
        decoration={markings, mark=at position -45/720 with {\arrow{>}}},
        postaction={decorate}
        ] (0,0) circle (\R cm);

\draw[
        decoration={markings, mark=at position (15+240)/720 with {\arrow{>}}},
        decoration={markings, mark=at position (45+240)/720 with {\arrow{>}}},
        decoration={markings, mark=at position (75+240)/720 with {\arrow{>}}},
        decoration={markings, mark=at position (-15+240)/720 with {\arrow{>}}},
        decoration={markings, mark=at position (-45+240)/720 with {\arrow{>}}},
        postaction={decorate}
        ] (0,0) circle (\R cm);

\draw[
        decoration={markings, mark=at position (15+480)/720 with {\arrow{>}}},
        decoration={markings, mark=at position (45+480)/720 with {\arrow{>}}},
        decoration={markings, mark=at position (75+480)/720 with {\arrow{>}}},
        postaction={decorate}
        ] (0,0) circle (\R cm);

\draw[
        decoration={markings, mark=at position 31/48 with {\arrow{>}}},
        postaction={decorate}
        ] (0,0) circle (\R cm);

\draw[
        decoration={markings, mark=at position 435/720 with {\arrow{>}}},
        postaction={decorate}
        ] (0,0) circle (\R cm);


\draw ({\R*cos(0/2)},{\R*sin(0/2)}) node {$\bullet$};
\draw ({\R*cos(30/2)},{\R*sin(30/2)}) node {$\bullet$};
\draw ({\R*cos(60/2)},{\R*sin(60/2)}) node {$\bullet$};
\draw ({\R*cos(-30/2)},{\R*sin(-30/2)}) node {$\bullet$};

\draw ({\R*(cos(15/2)+0.05)},{\R*(sin(15/2)-0)}) node[anchor=west] {$V_1$};
\draw ({\R*(cos(45/2)+0.05)},{\R*(sin(45/2)+0)}) node[anchor=west] {$V_2$};
\draw ({\R*(cos(75/2)-0)},{\R*(sin(75/2)+0.05)}) node[anchor=west] {\rotatebox{-52.5}{$\cdots$}};

\draw ({\R*(cos(-15/2)+0.05)},{\R*(sin(-15/2)-0)}) node[anchor=west] {$V_n$};
\draw ({\R*(cos(-45/2)-0)},{\R*(sin(-45/2)-0.05)}) node[anchor=west] {\rotatebox{67.5}{$\cdots$}};


\draw ({\R*cos(0/2+120)},{\R*sin(0/2+120)}) node {$\bullet$};
\draw ({\R*cos(30/2+120)},{\R*sin(30/2+120)}) node {$\bullet$};
\draw ({\R*cos(60/2+120)},{\R*sin(60/2+120)}) node {$\bullet$};
\draw ({\R*cos(-30/2+120)},{\R*sin(-30/2+120)}) node {$\bullet$};

\draw ({\R*(cos(15/2+120)-0.05)},{\R*(sin(15/2+120)+0.05)}) node[anchor=south] {$V_1$};
\draw ({\R*(cos(45/2+120)-0.1)},{\R*(sin(45/2+120)+0)}) node[anchor=south] {$V_2$};
\draw ({\R*(cos(75/2+120)-0)},{\R*(sin(75/2+120)+0.05)}) node[anchor=east] {\rotatebox{247.5}{$\cdots$}};

\draw ({\R*(cos(-15/2+120)+0)},{\R*(sin(-15/2+120)+0.05)}) node[anchor=south] {$V_n$};
\draw ({\R*(cos(-45/2+120)-0},{\R*(sin(-45/2+120)+0.05)}) node[anchor=south] {\rotatebox{7.5}{$\cdots$}};


\draw ({\R*cos(0/2+240)},{\R*sin(0/2+240)}) node {$\bullet$};
\draw ({\R*cos(30/2+240)},{\R*sin(30/2+240)}) node {$\bullet$};
\draw ({\R*cos(60/2+240)},{\R*sin(60/2+240)}) node {$\bullet$};
\draw ({\R*cos(-30/2+240)},{\R*sin(-30/2+240)}) node {$\bullet$};

\draw ({\R*(cos(15/2+240)-0.05)},{\R*(sin(15/2+240)-0.05)}) node[anchor=north] {$V_1$};
\draw ({\R*(cos(45/2+240)+0)},{\R*(sin(45/2+240)-0.05)}) node[anchor=north] {$V_2$};
\draw ({\R*(cos(75/2+240)-0)},{\R*(sin(75/2+240)-0.05)}) node[anchor=north] {\rotatebox{187.5}{$\cdots$}};

\draw ({\R*(cos(-15/2+240)-0.05)},{\R*(sin(-15/2+240)-0)}) node[anchor=north] {$V_n$};
\draw ({\R*(cos(-45/2+240)-0},{\R*(sin(-45/2+240)-0.05)}) node[anchor=east] {\rotatebox{127.5}{$\cdots$}};

\end{tikzpicture}
}
\right)^{\htpy \Cyclic_r} \right)
&
F
\arrow[mapsto]{l}
\end{tikzcd}
~,
\end{equation}  
which takes a functor $F$ to the new functor that precomposes with pullback along an $r$-fold cover and postcomposes with homotopy fixedpoints for the resulting $\Cyclic_r$-action.  In other words, this cartesian fibration (the downwards functor in diagram \Cref{global.section.of.relfun}) is classified by a right $\Nx$-action on the fiber \Cref{fiber.of.cartesian.fibn.over.BN}, through which the element $r \in \Nx$ acts as the endofunctor \Cref{cartesian.mdrmy.by.h}.  To emphasize that this right action involves \textit{homotopy} fixedpoints, we denote it by
\begin{equation}
\label{Nx.action.by.h}
\Fun \left( \int_{|\cyclic} \BV , \cV \right)
\ractsh
\Nx
~.
\end{equation}

\subsection{Cartesian monodromy via cocartesian monodromy in families}
\label{subsection.cart.mdrmy.via.cocart.mdrmy}

We now describe the cartesian monodromy functor \Cref{cartesian.mdrmy.by.h} more precisely.  For this, let us first recall that an object of the fiber
\[ \begin{tikzcd}
\int_{|\cyclic} \BV
\arrow[hook, two heads]{r}
\arrow{d}
&
\int_{|\epicyc} \BV
\arrow{d}
\\
\pt
\arrow[hook, two heads]{r}
&
\BN
\end{tikzcd} \]
is given by the data
\begin{equation}
\label{obj.of.int.epicyc.BV}
\left(
S \in \BW
\ , \
(R \ra S) \in \D(S)
\ , \
R^{(1)} \ra \cV
\right)
\end{equation}
of a circle, a disk-refinement thereof, and a labeling of the intervals of the latter by objects of $\cV$.  As both maps in the composite
\[ \begin{tikzcd}
\int_{|\epicyc} \BV
\arrow{d}
\\
\epicyc
\arrow{d}
\\
\BW
\end{tikzcd} \]
are cocartesian fibrations, then so is their composite.  Hence, for any map $S \ra S'$ in $\BW$ of degree $r \in \Nx$, we obtain a cocartesian pushforward
\[
\left(
S' \in \BW
\ , \
(R' \ra S') \in \D(S)
\ , \
R'^{(1)} \ra R^{(1)} \ra \cV
\right)
\]
of the object \Cref{obj.of.int.epicyc.BV}: the disk-refinement $R'$ is pulled back from $R$, and this induces an $r$-fold covering map $R'^{(1)} \ra R^{(1)}$ of discrete spaces (i.e.\! we pull back the labelings of intervals by objects of $\cV$ from $R$ to $R'$ as well).\footnote{In fact, we have already seen an example of such a cocartesian pushforward, namely the morphism \Cref{pic.of.morphism.in.int.BV.over.three.fold.creation}.}

Now, the space of degree-$r$ maps in $\BW$ with source $S$ is a copy of $\BC_r$.  Hence, cocartesian pushforward \textit{over all of these maps simultaneously} defines a $\Cyclic_r$-object
\[
\pi_r^* \Cref{obj.of.int.epicyc.BV}
:
\BC_r
\longra
\int_{|\cyclic} \BV
~.
\]
This construction assembles into a functor
\begin{equation}
\label{pi.r.notation.for.BV}
\pi_r^*
:
\int_{|\cyclic} \BV
\longra
\Fun \left(\BC_r,\int_{|\cyclic} \BV\right)
~,
\end{equation}
and with this notation we can rewrite the cartesian monodromy functor \Cref{cartesian.mdrmy.by.h} more formally as
\begin{equation}
\label{cartesian.mdrmy.by.h.revisited}
\begin{tikzcd}[row sep=0cm]
\Fun \left( \int_{|\cyclic} \BV , \cV \right)
&
\Fun \left( \int_{|\cyclic} \BV , \cV \right)
\arrow{l}
\\
\rotatebox{90}{$\in$}
&
\rotatebox{90}{$\in$}
\\
\left( F ( \pi_r^*(-) ) \right)^{\htpy \Cyclic_r}
&
F
\arrow[mapsto]{l}
\end{tikzcd}
~.
\end{equation}

This identification of cartesian monodromy in a relative functor $\infty$-category via ``cocartesian monodromy in families'' admits a useful generalization.  Namely, let
\[
( \cE \da \BW )
\in
\coCart_\BW
\]
be an arbitrary cocartesian fibration over $\BW$, and write 
\[
\begin{tikzcd}
\cE_0
\arrow[hook, two heads]{r}
\arrow{d}
&
\cE
\arrow{d}
\\
\BT
\arrow[hook, two heads]{r}
&
\BW
\end{tikzcd}
\]
for the pullback (or equivalently, the fiber over the unique point in $\BN$).  Let us repurpose our notation $\pi_r^*$ for the ``all the cocartesian monodromy over degree-$r$ maps in $\BW$ at once'' functor \Cref{pi.r.notation.for.BV} to still write
\[
\pi_r^*
:
\cE_0
\longra
\Fun(\BC_r,\cE_0)
\]
in this more general setting.  Then, by [\Cref{cyclo:identify.relFun.as.cart.fibn.with.htpy.fixedpts}], the functor
\begin{equation}
\label{Fun.rel.BN.E.V}
\begin{tikzcd}
\Fun^\rel_{/\BN}(\cE,\ul{\cV})
\arrow{d}
\\
\BN
\end{tikzcd}
\end{equation}
is once again a cartesian fibration: its fiber is $\Fun(\cE_0,\cV)$, and the cartesian monodromy functor over the morphism $[1] \xra{r} \BN$ is given by the formula
\[
\begin{tikzcd}[row sep=0cm]
\Fun(\cE_0,\cV)
&
\Fun(\cE_0,\cV)
\arrow{l}
\\
\rotatebox{90}{$\in$}
&
\rotatebox{90}{$\in$}
\\
(F(\pi_r^*(-)))^{\htpy \Cyclic_r}
&
F
\arrow[mapsto]{l}
\end{tikzcd}
~.
\]
We denote this right $\Nx$-module by
\begin{equation}
\label{general.notation.for.right.Nx.action.on.Fun.Ezero.V}
\Fun(\cE_0,\cV)
\ractsh
\Nx
~,
\end{equation}
extending the notation \Cref{Nx.action.by.h}.  Evidently, this construction assembles into a functor
\begin{equation}
\label{relFun.BW.to.Nx}
\Fun^\rel_{/\BN}(-,\ul{\cV})
:
\left( \coCart_\BW \right)^\op
\longra
\Cart_\BN
~.
\end{equation}

\subsection{The diagonal package}
\label{subsection.diagonal.package}

Let us now return to the global section \Cref{global.section.of.relfun}.  By definition, this is an object
\begin{equation}
\label{diagonal.package}
\left(
\prod
,
\htpy
\right)
\in
\lim^\rlax
\left(
\Fun \left( \int_{|\cyclic} \BV , \cV \right)
\ractsh
\Nx
\right)
:=
\Gamma
\left(
\begin{tikzcd}
\Fun^\rel_{/\BN} \left( \int_{|\epicyc} \BV , \ul{\cV} \right)
\arrow{d}
\\
\BN
\end{tikzcd}
\right)
\simeq
\Fun
\left(
\int_{|\epicyc} \BV
,
\cV
\right)
\end{equation}
of the right-lax limit of the right $\Nx$-module \Cref{Nx.action.by.h}: over the unique object of $\BN$ this selects the object
\[
\prod
\in
\Fun
\left(
\int_{|\cyclic} \BV,
\cV
\right)
~,
\]
and over a morphism $[1] \xra{r} \BN$ this selects the diagonal map
\[
\prod
\longra
\left( \prod \pi_r^*(-)\right)^{\htpy \Cyclic_r}
~.\footnote{In fact, this map is an equivalence: for any objects $V_1,\ldots,V_k \in \cV$, the map
\[
\prod_i V_i
\longra
\left( \prod_i (V_i)^{\times r} \right)^{\htpy \Cyclic_r}
\]
is an equivalence.  In other words, the global section \Cref{global.section.of.relfun} actually lands in the subcategory of \textit{cartesian} morphisms, and hence defines a point in the \textit{strict} limit: the $\infty$-category of cartesian sections of the corresponding cartesian fibration.  However, this will not be relevant for us.}
\]
We include the ``$\htpy$'' in the notation of the object \Cref{diagonal.package} to emphasize that these structure maps take place within the ambient context of a right $\Nx$-action involving homotopy fixedpoints.  We refer to the object \Cref{diagonal.package} as the \bit{diagonal package}, since as we will see in \Cref{subsection.unstable.cyclo.str} it encodes diagonal maps in $\cV$ along with all of their requisite compatibilities necessary for defining the $\WW$-action \Cref{fact.hlgy.with.W.action} on $\THHVx$.

\subsection{Unstable cyclotomic objects in $\cV$}
\label{subsection.unstable.cyclo.objects}

Consider the functor
\[
\begin{tikzcd}
\Fun^\rel_{/\BN}(\BW,\ul{\cV})
\arrow{d}
\\
\BN
\end{tikzcd}
~.
\]
By the discussion of \Cref{subsection.cart.mdrmy.via.cocart.mdrmy}, this is a cartesian fibration, whose fiber over the unique object of $\BN$ is
\[
\Fun ( \BT,\cV )
\]
and whose cartesian monodromy over the morphism $[1] \xra{r} \BN$ is
\[ \begin{tikzcd}[row sep=0cm]
\Fun(\BT,\cV)
&
\Fun(\BT,\cV)
\arrow{l}
\\
\rotatebox{90}{$\in$}
&
\rotatebox{90}{$\in$}
\\
\left(
\TT
\simeq
(\TT / \Cyclic_r)
\lacts
V^{\htpy \Cyclic_r}
\right)
&
\left(
\TT
\lacts
V
\right)
\arrow[mapsto]{l}
\end{tikzcd}~; \]
according to our notational convention \Cref{general.notation.for.right.Nx.action.on.Fun.Ezero.V}, we write
\begin{equation}
\label{Nx.mod.whose.rlax.lim.is.unstable.cyclo.objs}
\Fun(\BT,\cV)
\ractsh
\Nx
\end{equation}
for the corresponding right $\Nx$-module.  Then, we define the $\infty$-category of \bit{unstable cyclotomic objects} in $\cV$ to be
\[
\Cyclo^\htpy(\cV)
:=
\lim^\rlax
\left(
\Fun(\BT,\cV)
\ractsh
\Nx
\right)
:=
\Gamma
\left(
\begin{tikzcd}
\Fun^\rel_{/\BN}(\BW,\ul{\cV})
\arrow{d}
\\
\BN
\end{tikzcd}
\right)
\simeq
\Fun(\BW,\cV)
~,
\]
the right-lax limit of the right $\Nx$-module \Cref{Nx.mod.whose.rlax.lim.is.unstable.cyclo.objs}.\footnote{One might alternatively refer to this as the $\infty$-category of \bit{cyclohomic objects}, since an unstable cyclotomic structure involves \textit{homotopy} fixedpoints whereas a \textit{cyclotomic} structure on a spectrum will involve \textit{Tate} fixedpoints.}\footnote{The simplest example of an unstable cyclotomic object is a \textit{free loop object}, i.e.\! the cotensoring
\[ S^1 \cotensoring V \]
of the circle into an arbitrary object $V \in \cV$ (e.g.\! the free loopspace of a space); the unstable cyclotomic structure arises from the identifications $S^1 \simeq (S^1)_{\htpy \Cyclic_r}$ (compatibly for all $r \in \Nx$), so that its $\TT$-equivariant structure maps
\[
S^1 \cotensoring V
\longra
(S^1 \cotensoring V)^{\htpy \Cyclic_r}
\simeq (S^1)_{\htpy \Cyclic_r} \cotensoring V
\]
are actually \textit{equivalences}.  When $\cV = \Spaces$, this is nothing but $\THH^\times_\Spaces(V)$, the spatially-enriched factorization homology of the $\infty$-groupoid $V \in \Spaces \subset \Cat$ over the circle.}  In analogy with our notation for cyclotomic spectra, we denote an unstable cyclotomic structure on an object $V \in \cV$ by
\[
\Cyclo^\htpy
\lacts
V
~.
\]

\subsection{The unstable cyclotomic structure on $\THHVx$}
\label{subsection.unstable.cyclo.str}

We now explain how the diagonal package \Cref{diagonal.package} gives rise to the unstable cyclotomic structure on $\THHVx(\cC)$.

We begin by noting that the diagram
\[ \begin{tikzcd}[column sep=2cm, row sep=1.5cm]
\int_{|\epicyc} \cdiC
\arrow{r}{\int_{|\epicyc} \ulhom_\cC}
\arrow{d}
&
\int_{|\epicyc} \BV
\arrow{ld}
\\
\epicyc
\arrow{d}
\\
\BW
\end{tikzcd} \]
defines a morphism in $\coCart_\BW$ [\Cref{fact:obs.cocart.over.cls.and.some.cr}].  Trivially, we obtain a span
\begin{equation}
\label{span.in.cocart.BW}
\begin{tikzcd}
\int_{|\epicyc} \BV
\arrow{rd}
&
\int_{|\epicyc} \cdiC
\arrow{l}
\arrow{r}
\arrow{d}
&
\BW
\arrow{ld}[sloped, pos=0.6]{\sim}
\\
&
\BW
\end{tikzcd}
\end{equation}
in $\coCart_\BW$.  Applying the functor \Cref{relFun.BW.to.Nx} to this span gives a cospan
\begin{equation}
\label{cospan.in.Cart.BN}
\begin{tikzcd}
\Fun^\rel_{/\BN} \left( \int_{|\epicyc} \BV , \ul{\cV} \right)
\arrow{r}
\arrow{rd}
&
\Fun^\rel_{/\BN} \left( \int_{|\epicyc} \cdiC , \ul{\cV} \right)
\arrow{d}
&
\Fun^\rel_{/\BN} ( \BW , \ul{\cV} )
\arrow{l}
\arrow{ld}
\\
&
\BN
\end{tikzcd}
\end{equation}
in $\Cart_\BN$.

Now, on fibers, the leftward functor of diagram \Cref{cospan.in.Cart.BN} restricts to a functor admitting a left adjoint
\[ \begin{tikzcd}[column sep=1.5cm]
\Fun \left( \int_{|\cyclic} \cdiC , \cV \right)
\arrow[dashed, transform canvas={yshift=0.9ex}]{r}
\arrow[leftarrow, transform canvas={yshift=-0.9ex}]{r}[transform canvas={yshift=-0.2ex}]{\bot}
&
\Fun ( \BW , \cV )
\end{tikzcd}
~,
\]
namely the left Kan extension appearing in diagram \Cref{fact.hlgy.with.T.action}.  By completely general principles [\Cref{cyclo:lem.get.r.lax.left.adjt}], it follows that there exists a (necessarily unique) extension
\begin{equation}
\label{extend.ladjt.to.radjt.in.Cat.cart.over.BN}
\begin{tikzcd}[column sep=1.5cm, row sep=1.5cm]
\Fun \left( \int_{|\cyclic} \cdiC , \cV \right)
\arrow[transform canvas={yshift=0.9ex}]{r}
\arrow[leftarrow, transform canvas={yshift=-0.9ex}]{r}[transform canvas={yshift=-0.2ex}]{\bot}
\arrow[hook]{d}
&
\Fun ( \BW , \cV )
\arrow[hook]{d}
\\
\Fun^\rel_{/\BN} \left( \int_{|\epicyc} \cdiC , \ul{\cV} \right)
\arrow[dashed, transform canvas={yshift=0.9ex}]{r}
\arrow[leftarrow, transform canvas={yshift=-0.9ex}]{r}[transform canvas={yshift=-0.2ex}]{\bot}
&
\Fun^\rel_{/\BN} ( \BW , \ul{\cV} )
\end{tikzcd}
\end{equation}
of this left adjoint, a morphism in $\Cat_{\cart/\BN}$.  In fact, the lower row of diagram \Cref{extend.ladjt.to.radjt.in.Cat.cart.over.BN} defines an adjunction in the $(\infty,2)$-category $\Cat_{\cart/\BN}$ (which justifies the notation), i.e.\! a functor
\begin{equation}
\label{adjn.in.Cat.cart.over.BN}
\Adj
\longra
\Cat_{\cart/\BN}
\end{equation}
of $(\infty,2)$-categories from the walking adjunction (which is actually just a $2$-category).  Then, postcomposing the functor \Cref{adjn.in.Cat.cart.over.BN} with the global sections functor
\[
\Cat_{\cart/\BN}
\xlongra{\Gamma}
\Cat
\]
of $(\infty,2)$-categories yields an adjunction
\[ \begin{tikzcd}[column sep=1.5cm]
\Fun \left( \int_{|\epicyc} \cdiC , \cV \right)
\arrow[transform canvas={yshift=0.9ex}]{r}
\arrow[leftarrow, transform canvas={yshift=-0.9ex}]{r}[transform canvas={yshift=-0.2ex}]{\bot}
&
\Fun ( \BW , \cV )
\end{tikzcd} \]
in $\Cat$, whose left adjoint is the left Kan extension appearing in diagram \Cref{fact.hlgy.with.W.action}.  Hence, we obtain the unstable cyclotomic structure on $\THHVx(\cC)$ as the image of the diagonal package, as illustrated in \Cref{diagonal.package.gives.W.action}.
\begin{sidewaysfigure}
\vspace{425pt}
\[ \begin{tikzcd}[row sep=2.5cm, column sep=2cm]
\Fun \left( \int_{|\cyclic} \BV, \cV \right)
\arrow{r}
\arrow[leftarrow]{d}
&
\Fun \left( \int_{|\cyclic} \cdiC , \cV \right)
\arrow[dashed, transform canvas={yshift=0.9ex}]{r}
\arrow[leftarrow, transform canvas={yshift=-0.9ex}]{r}[transform canvas={yshift=-0.2ex}]{\bot}
\arrow[leftarrow]{d}
&
\Fun ( \BT , \cV )
\arrow[leftarrow]{d}
\\
\Fun \left( \int_{|\epicyc} \BV, \cV \right)
\arrow{r}
&
\Fun \left( \int_{|\epicyc} \cdiC , \cV \right)
\arrow[dashed, transform canvas={yshift=0.9ex}]{r}
\arrow[leftarrow, transform canvas={yshift=-0.9ex}]{r}[transform canvas={yshift=-0.2ex}]{\bot}
&
\Fun ( \BW , \cV )
\\[-2.5cm]
\rotatebox{90}{$\simeq$}
&
\rotatebox{90}{$\simeq$}
&
\rotatebox{90}{$\simeq$}
\\[-2.5cm]
\lim^\rlax
\left(
\Fun \left( \int_{|\cyclic} \BV , \cV \right)
\ractsh
\Nx
\right)
\arrow{r}
&
\lim^\rlax
\left(
\Fun \left( \int_{|\cyclic} \cdiC , \cV \right)
\ractsh
\Nx
\right)
\arrow[dashed, transform canvas={yshift=0.9ex}]{r}
\arrow[leftarrow, transform canvas={yshift=-0.9ex}]{r}[transform canvas={yshift=-0.2ex}]{\bot}
&
\lim^\rlax
\left(
\Fun \left( \BT , \cV \right)
\ractsh
\Nx
\right)
&
\hspace{-2.5cm} =: \Cyclo^\htpy(\cV)
\\
\prod
\arrow[maps to]{rr}{\Cref{fact.hlgy.with.T.action}}
&
&
\left( \TT \lacts \THHVx(\cC) \right)
\\
\prod
\arrow[maps to]{u}
\arrow[maps to]{rr}{\Cref{fact.hlgy.with.W.action}}
&
&
\left( \WW \lacts \THHVx(\cC) \right)
\arrow[mapsto]{u}
\\[-1.75cm]
\left( \prod , \htpy \right)
\arrow[leftrightarrow]{u}
\arrow[maps to]{rr}
&
&
\left( \Cyclo^\htpy \lacts \THHVx(\cC) \right)
\arrow[leftrightarrow]{u}
\end{tikzcd} \]
\vspace{2.5cm}
\caption{Above: the diagram of $\infty$-categories housing the passage from the diagonal package to the unstable cyclotomic structure on $\THHVx(\cC)$.  All solid horizontal morphisms are given by pullback; the upwards functors are given by restriction to underlying objects of right $\Nx$-modules, and the diagram commutes when omitting either all left adjoints or all right adjoints.  Below: the diagonal package gives rise to the unstable cyclotomic structure on $\THHVx(\cC)$, and the restriction of this passage recovers its underlying $\TT$-action.}
\label{diagonal.package.gives.W.action}
\end{sidewaysfigure}

In \Cref{unstable.cyclo.str.via.modules}, we summarize the constructions of this subsection in the language of modules of [\Cref{cyclo:section.lax.actions.and.limits}]; our method for obtaining the cyclotomic structure on $\THH$ will be closely analogous.
\begin{figure}[h]
\begin{tikzcd}[column sep=4.5cm, row sep=1.5cm]
{[1] \underset{1,\pt,1}{\coprod} [1]}
\arrow{r}{\left( \int_{|\epicyc} \BV \la \int_{|\epicyc} \cdiC \ra \BW \right)^\circ}
\arrow{d}[swap]{\id \underset{\pt}{\coprod} {\sf r.adjt}}  
&
\left( \LMod_\WW \right)^\op
\arrow{r}{\left( \Fun ( (-)_0 , \cV) \ractsh \Nx \right)}
&
\RMod_\Nx
\arrow[hook, two heads]{d}
\\
{ [1] \underset{1,\pt,l}{\coprod} \Adj}
\arrow[dashed]{rr}{\exists ! ~ [\textup{\Cref{cyclo:lem.get.r.lax.left.adjt}}]}
\arrow[leftarrow]{d}[swap]{\id \underset{\pt}{\coprod} {\sf l.adjt}} 
&
&
\RMod^\rlax_\Nx
\arrow{dd}{\lim^\rlax}
\\
{[2]}
\arrow[leftarrow]{d}[swap]{\{0 < 2 \}}  
\\
{[1]}
\arrow[dashed]{rr}
&
&
\Cat
\end{tikzcd}
\caption{A diagram illustrating in the language of modules how to obtain the lower rightwards composite functor in \Cref{diagonal.package.gives.W.action}, which takes the diagonal package to the unstable cyclotomic structure on $\THHVx(\cC)$.  Here, we write $l \in \Adj$ for the source of the walking left adjoint.}
\label{unstable.cyclo.str.via.modules}
\end{figure}

\subsection{From fiberwise colimits to unstable cyclotomic structures}
\label{subsection.unwind.unst.cyclo.str.map}

Purely for expository purposes, we unwind here how the fiberwise colimit of diagram \Cref{fact.hlgy.with.W.action} gives rise to the unstable cyclotomic structure maps
\begin{equation}
\label{unst.cyclo.str.map.for.cC}
\THHVx(\cC)
\longra
\THHVx(\cC)^{\htpy \Cyclic_r}
~;
\end{equation}
this is essentially just an exercise in the functoriality of fiberwise colimits along cocartesian fibrations.  For this, let 
\[
S_s
\xlongra{\varphi}
S_t
\]
be a morphism in $\BW$ corresponding contravariantly to an $r$-fold covering map of framed circles, and write 
\[
\D(S_s)
\xlongra{\varphi_*}
\D(S_t)
\]
for the resulting functor on categories of disk-refinements (the \textit{$r$-fold subdivision} functor on $\para$).  Then, the map \Cref{unst.cyclo.str.map.for.cC} is given by the composite
\begin{align}
\nonumber
\THHVx(\cC)
& :=
\int_{S_s} \cC
\\
\nonumber
& :=
\underset{\substack{(R_s \ra S_s) \in \D(S_s) \\ (C_i) \in (\iC)^{\times R_s^{(0)}} }}{\colim} \left( \prod_{i \in \ZZ/R_s^{(0)}} \ulhom_\cC(C_i , C_{i+1}) \right)
\\
\label{diagonal.map.for.unst.cyclo.str}
& \longra
\underset{\substack{(R_s \ra S_s) \in \D(S_s) \\ (C_i) \in (\iC)^{\times R_s^{(0)}} }}{\colim} \left( \left( \prod_{i \in \ZZ/R_s^{(0)}} \ulhom_\cC(C_i , C_{i+1})^{\times r} \right)^{\htpy \Cyclic_r} \right)
\\
\label{interchange.map.for.unst.cyclo.str}
& \longra
\left( \underset{\substack{(R_s \ra S_s) \in \D(S_s) \\ (C_i) \in (\iC)^{\times R_s^{(0)}} }}{\colim} \left( \prod_{i \in \ZZ/R_s^{(0)}} \ulhom_\cC(C_i , C_{i+1})^{\times r} \right) \right)^{\htpy \Cyclic_r}
\\
\label{inclusion.map.for.unst.cyclo.str}
& \longra
\left( \underset{\substack{(R_t \ra S_t) \in \D(S_t) \\ (C_j) \in (\iC)^{\times R_t^{(0)}} }}{\colim} \left( \prod_{j \in \ZZ/R_t^{(0)}} \ulhom_\cC(C_j , C_{j+1}) \right) \right)^{\htpy \Cyclic_r}
\\
\nonumber
& =:
\left( \int_{S_t} \cC \right)^{\htpy \Cyclic_r}
\\
\nonumber
& =:
\left( \THHVx(\cC) \right)^{\htpy \Cyclic_r}
~,
\end{align}
in which
\begin{itemize}
\item the map \Cref{diagonal.map.for.unst.cyclo.str} arises from the diagonal maps in $\cV$ (and is actually an equivalence),
\item the map \Cref{interchange.map.for.unst.cyclo.str} arises from the universal property of colimits, and
\item the map \Cref{inclusion.map.for.unst.cyclo.str} arises from a map on colimits taking the term indexed by
\[
\left(
R_s \in \D(S_s)
\ , \
R_s^{(0)} \xra{(C_i)} \iC
\right)
\]
to the equivalent term indexed by
\[
\left(
\varphi_*(R_s) \in \D(S_t)
\ , \
(\varphi_*(R_s))^{(0)}
\longra
R_s^{(0)}
\xra{(C_i)}
\iC
\right)
~.
\]
\end{itemize}

\section{The cyclotomic structure on $\THH$}
\label{section.roadmap.cyclo.str}

In this section, we specialize to the case that $(\cV,\boxtimes) = (\Spectra , \otimes)$ is the $\infty$-category of spectra equipped with its smash product monoidal structure, and we obtain the cyclotomic structure
\[
\Cyclo
\lacts
\THH
~.
\]
We achieve this in \Cref{subsection.cyclo.str.on.THH}, in a nearly identical way to how we obtained the unstable cyclotomic structure on $\THHVx$ in \Cref{subsection.unstable.cyclo.str}.  In particular, whereas the unstable cyclotomic structure on $\THHVx$ arose from the \textit{diagonal package}, the cyclotomic structure on $\THH$ will arise from the \bit{Tate package}; we explain this analogy further in \Cref{subsection.locate.tate.package}.  We will obtain the Tate package from the diagonal package via linearization (in the sense of Goodwillie calculus); we achieve this in \Cref{subsection.tate.package}, after describing certain auxiliary constructions that collectively accommodate this maneuver in \Cref{subsection.aux.constrns}.

\subsection{Locating the Tate package}
\label{subsection.locate.tate.package}

As described in \Cref{subsection.unwind.unst.cyclo.str.map}, the unstable cyclotomic structure on $\THHVx$ consists of a system of maps
\[
\THHVx
\longra
( \THHVx )^{\htpy \Cyclic_r}
\]
which are ultimately induced by the natural diagonal maps
\[ \begin{tikzcd}
\prod_i V_i
\arrow{rr}
\arrow[dashed]{rd}[swap, sloped, pos=0.6]{\sim}
&
&
\prod_i (V_i)^{\times r}
\\
&
\left( \prod_i (V_i)^{\times r} \right)^{\htpy \Cyclic_r}
\arrow{ru}
\end{tikzcd}
\]
that exist for any objects $V_1,\ldots, V_k \in \cV$.  In \Cref{subsection.diagonal.package}, we organized these maps into the \textit{diagonal package}, an object
\[
\left( \prod , \htpy \right)
\in
\lim
\left(
\Fun \left( \int_{|\cyclic} \BV , \cV \right)
\ractsh
\Nx
\right)
\]
(originally \Cref{diagonal.package}) in the limit of a certain right $\Nx$-action by homotopy fixedpoints in which the element $r \in \Nx$ acts as the endofunctor
\[ \begin{tikzcd}[row sep=0cm]
\Fun \left( \int_{|\cyclic} \BV , \cV \right)
&
\Fun \left( \int_{|\cyclic} \BV , \cV \right)
\arrow{l}
\\
\rotatebox{90}{$\in$}
&
\rotatebox{90}{$\in$}
\\
\left( F ( \pi_r^*(-) ) \right)^{\htpy \Cyclic_r}
&
F
\arrow[mapsto]{l}
\end{tikzcd} \]
(originally \Cref{cartesian.mdrmy.by.h.revisited}).

By contrast, the cyclotomic structure on $\THH$ will be induced in a formally identical way by the natural \textit{Tate diagonal} maps
\[ \begin{tikzcd}
\bigotimes_i E_i
\arrow[dotted]{rr}[description]{\not\exists}
\arrow[dotted]{rd}[description]{\not\exists}
\arrow[dashed]{rdd}
&
&
\bigotimes_i (E_i)^{\otimes r}
\\
&
\left( \bigotimes_i (E_i)^{\otimes r} \right)^{\htpy \Cyclic_r}
\arrow{ru}
\arrow{d}
\\
&
\left( \bigotimes_i (E_i)^{\otimes r} \right)^{\tate \Cyclic_r}
\end{tikzcd} \]
that exist for any spectra $E_1,\ldots,E_k \in \Spectra$.  
But in fact, it turns out that the formula
\[ \begin{tikzcd}[row sep=0cm]
\Fun \left( \int_{|\cyclic} \BSp , \Spectra \right)
&
\Fun \left( \int_{|\cyclic} \BSp , \Spectra \right)
\arrow{l}
\\
\rotatebox{90}{$\in$}
&
\rotatebox{90}{$\in$}
\\
\left( F ( \pi_r^*(-) ) \right)^{\tate \Cyclic_r}
&
F
\arrow[mapsto]{l}
\end{tikzcd} \]
only defines a \textit{left-lax} right $\Nx$-action.  Thus we find that the main ingredient in obtaining the cyclotomic structure on $\THH$ will be the \bit{Tate package}, an object
\[
\left( \bigotimes , \tate \right)
\in
\lim^\rlax
\left(
\Fun \left( \int_{|\cyclic} \BSp , \Spectra \right)
\ractstau
\Nx
\right)
\]
of the right-lax limit of a certain left-lax right $\Nx$-action by Tate fixedpoints.

\subsection{Auxiliary constructions}
\label{subsection.aux.constrns}

In order to obtain the Tate package, we make the following three constructions; the first two follow from [\Cref{cyclo:proto.tate}],\footnote{Using the language of \cite{AMR-cyclo}, these two deductions are supported by the functor
\[
\coCart_\BW
\longra
\Fun \left( \BNop ,\GSpan \right)
\]
which takes $(\cE \da \BW) \in \coCart_\BW$ to the functor $\BNop \ra \GSpan$ defined as follows:
\begin{itemize}
\item it takes the unique object of $\BNop$ to the object $\cE_0 \in \GSpan$, and
\item it takes the morphism $[1] \xra{r^\circ} \BNop$ to the span
\[ \begin{tikzcd}[ampersand replacement=\&]
\&
\Gamma_{[1]}^{\cocart/\BW}(\cE)
\arrow{ld}
\arrow{rd}[pos=0.4]{\BC_r\textup{-}{\sf Kan}}
\\
\Gamma_{\{0\}}(\cE)
\&
\&
\Gamma_{\{1\}}(\cE)
\end{tikzcd} \]
from the full subcategory $\Gamma_{[1]}^{\cocart/\BW}(\cE) \subset \Gamma_{[1]}(\cE)$ on those sections
\[ \begin{tikzcd}[ampersand replacement=\&]
\&
\cE
\arrow{d}
\\
\&
\BW
\arrow{d}
\\
{[1]}
\arrow{r}[swap]{r}
\arrow[dashed]{ruu}
\&
\BN
\end{tikzcd} \]
which select morphisms in $\cE$ that are cocartesian over their images in $\BW$.
\end{itemize}
(Note that this prescription respects composition since $\cE \da \BW$ is a cocartesian fibration, i.e.\! the cocartesian morphisms are stable under composition.)} while the third is simply a definition.  As in \Cref{subsection.cart.mdrmy.via.cocart.mdrmy}, we write $(\cE \da \BW) \in \coCart_\BW$ for an arbitrary object, and we write $\cE_0 = \cE_{|\BT}$.
\begin{enumerate}

\item We construct a left-lax right $\Nx$-action on
\[ \Fun \left( \cE_0 , \Spectra \right) ~ , \]
under which the element $r \in \Nx$ acts by the functor
\[ \begin{tikzcd}[row sep=0cm]
\Fun \left( \cE_0 , \Spectra \right)
&
\Fun \left( \cE_0 , \Spectra \right)
\arrow{l}
\\
\rotatebox{90}{$\in$}
&
\rotatebox{90}{$\in$}
\\
\left( F ( \pi_r^*(-) ) \right)^{\tate \Cyclic_r}
&
F
\arrow[mapsto]{l}
\end{tikzcd}
~ . \]
We denote this by
\begin{equation}
\label{llax.action.by.tau}
\Fun \left( \cE_0 , \Spectra \right)
\ractstau
\Nx
~.
\end{equation}
This defines a functor
\begin{equation}
\label{cocart.over.BW.to.llax.r.Nx.mod}
\left( \coCart_\BW \right)^\op
\xra{\left( \Fun ( (-)_0 , \Spectra ) \ractstau \Nx \right)}
\RMod_{\llax.\Nx}
~.
\end{equation}

\item We construct a right-lax $\Nx$-equivariant functor
\begin{equation}
\label{h.to.tau}
\left(
\begin{tikzcd}
\Fun^\rel_{/\BN} ( \cE , \Spectra )
\arrow{d}
\\
\BN
\end{tikzcd}
\right)
=:
\left(
\Fun \left( \cE_0 , \Spectra \right)
\ractsh
\Nx
\right)
\xra{\rlax}
\left(
\Fun \left( \cE_0 , \Spectra \right)
\ractstau
\Nx
\right)
\end{equation}
of left-lax right $\Nx$-modules.  This is the identity functor
\[
\Fun(\cE_0,\Spectra)
\xlongra{\id}
\Fun(\cE_0,\Spectra)
\]
on their common underlying $\infty$-category $\Fun(\cE_0,\Spectra)$, and for an object
\[
F
\in
\Fun(\cE_0,\Spectra)
\]
(considered as living in the source), the right-laxness determines structure maps
\[
\id \left( F ( \pi_r^*(-) )^{\htpy \Cyclic_r} \right)
\simeq
F ( \pi_r^*(-) )^{\htpy \Cyclic_r}
\longra
F ( \pi_r^*(-) )^{\tate \Cyclic_r}
\simeq
\id(F) ( \pi_r^*(-) )^{\tate \Cyclic_r}
\]
in $\Spectra$ (where we use ``$\id$'' in our notation only to clarify where these various constructions are taking place).

\item For any symmetric monoidal $\infty$-category $\cV$ admitting pushouts, we define a full subcategory
\begin{equation}
\label{subcat.of.fiberwise.multilinear.functors}
\Lin \left( \int_{|\cyclic} \BV , \Spectra \right)
\subset
\Fun \left( \int_{|\cyclic} \BV , \Spectra \right)
\end{equation}
of \bit{fiberwise multilinear} functors, as follows.  Recall first that a functor is \textit{linear} if it takes pushouts to pullbacks, and a multifunctor is \textit{multilinear} if it is linear separately in each variable.  Then, we define the full subcategory \Cref{subcat.of.fiberwise.multilinear.functors} to consist of those functors
\[
\int_{|\cyclic} \BV
\xlongra{F}
\Spectra
\]
that become multilinear when restricted to every fiber
\[ \begin{tikzcd}
\cV^{\times R^{(1)}}
\arrow{r}
\arrow{d}
&
\int_{|\cyclic} \BV
\arrow{r}{F}
\arrow{d}
&
\Spectra
\\
\{ (R \ra S) \}
\arrow{r}
&
\cyclic
\end{tikzcd}
~.
\]

The subcategory \Cref{subcat.of.fiberwise.multilinear.functors} is stable under the left-lax right $\Nx$-action \Cref{llax.action.by.tau}: this is immediate from the fundamental calculus-theoretic feature of the Tate construction [\Cref{cyclo:Cr.genzd.tate.of.rth.power.is.exact}], that the endofunctor
\[
E
\longmapsto
\left( E^{\otimes r} \right)^{\tate \Cyclic_r}
\]
of $\Spectra$ is linear.  The construction of the subobject
\[
\left(
\Lin \left( \int_{|\cyclic} \BV , \Spectra \right)
\ractstau
\Nx
\right)
\subset
\left(
\Fun \left( \int_{|\cyclic} \BV , \Spectra \right)
\ractstau
\Nx
\right)
\]
is evidently contravariantly functorial in $\cV$ for linear symmetric monoidal functors.

\end{enumerate}

\subsection{The Tate package}
\label{subsection.tate.package}

In order to construct the Tate package, we will first construct the \textit{finite} Tate package, i.e.\! its restriction to finite spectra, and then we will extend this to the Tate package itself.  This two-step process is necessary because the endofunctor
\[
E
\longmapsto
\left( E^{\otimes r} \right)^{\tate \Cyclic_r}
\]
on $\Spectra$ is linear but does not preserve filtered colimits [\Cref{cyclo:gen.tate.not.filt.cocts}].

We construct the finite Tate package from the finite diagonal package (which is simply the restriction of the diagonal package to finite spaces) as in \Cref{tate.package.from.diagonal.package}.
\begin{sidewaysfigure}
\vspace{400pt}
\[
\hspace{-70pt}
\begin{tikzcd}[column sep=2.5cm, row sep=2.5cm]
\left(
\Fun \left( \int_{|\cyclic} \BS^\fin , \Spaces \right)
\ractsh
\Nx
\right)
\arrow{r}{\Sigma^\infty_+ \circ -}[swap]{\rlax}
&
\left(
\Fun \left( \int_{|\cyclic} \BS^\fin , \Spectra \right)
\ractsh
\Nx
\right)
\arrow{r}{\Cref{h.to.tau}}[swap]{\rlax}
&
\left(
\Fun \left( \int_{|\cyclic} \BS^\fin , \Spectra \right)
\ractstau
\Nx
\right)
&
\left(
\Fun \left( \int_{|\cyclic} \BSp^\fin , \Spectra \right)
\ractstau
\Nx
\right)
\arrow{l}[swap]{- \circ \left( \int_{|\cyclic} \fB \Sigma^\infty_+\right)}
\\
&
&
\left(
\Lin \left( \int_{|\cyclic} \BS^\fin , \Spectra \right)
\ractstau
\Nx
\right)
\arrow[hook]{u}
&
\left(
\Lin \left( \int_{|\cyclic} \BSp^\fin , \Spectra \right)
\ractstau
\Nx
\right)
\arrow[hook]{u}
\arrow{l}{\sim}
\\
\left( \prod , \htpy \right)^\fin
\arrow[mapsto]{r}
&
\left( \bigotimes \circ \int_{|\cyclic} \fB \Sigma^\infty_+ , \htpy \right)^\fin
\arrow[mapsto]{rd}
\\
&
&
\left( \bigotimes \circ \int_{|\cyclic} \fB \Sigma^\infty_+ , \uptau \right)^\fin
\arrow[mapsto]{r}
&
\left( \bigotimes , \uptau \right)^\fin
\end{tikzcd} \]
\vspace{80pt}
\caption{Above: a commutative diagram of right-lax equivariant functors among left-lax right $\Nx$-modules.  Below: upon passage to right-lax limits, the finite diagonal package $(\prod,\htpy)^\fin$ is sent to the finite Tate package $(\bigotimes,\uptau)^\fin$.
\label{tate.package.from.diagonal.package}
}
\end{sidewaysfigure}
A number of comments are in order.
\begin{enumerate}

\item Perhaps unexpectedly, it is possible to take the left- or right-lax limit of any left- or right-lax left or right $\Nx$-module (with all three instances of ``or'' interpreted independently); on the other hand, a right-lax $\Nx$-equivariant morphism only induces a map on right-lax limits (and dually), while a strict morphism induces maps on all limits [\Cref{cyclo:section.lax.actions.and.limits}].  Thus, it is indeed meaningful to apply the functor $\lim^\rlax$ to the diagram in \Cref{tate.package.from.diagonal.package}.

\item The first morphism
\[
\left(
\Fun \left( \int_{|\cyclic} \BS^\fin , \Spaces \right)
\ractsh
\Nx
\right)
\xra[\rlax]{\Sigma^\infty_+ \circ -}
\left(
\Fun \left( \int_{|\cyclic} \BS^\fin , \Spectra \right)
\ractsh
\Nx
\right)
\]
in \Cref{tate.package.from.diagonal.package} is induced by the covariant functoriality of the relative functor $\infty$-category in the second variable.  It is only right-lax $\Nx$-equivariant because the suspension spectrum functor $\Sigma^\infty_+$ doesn't commute with the homotopy fixedpoints functors $(-)^{\htpy \Cyclic_r}$: given an object
\[
F
\in
\Fun \left( \int_{|\cyclic} \BS^\fin , \Spaces \right)
\]
of the source, the right-laxness of this morphism determines structure maps
\[
\Sigma^\infty_+ \left( \left( F ( \pi_r^*(-) ) \right)^{\htpy \Cyclic_r} \right)
\longra
\left( \Sigma^\infty_+ \left( F ( \pi_r^*(-) ) \right) \right)^{\htpy \Cyclic_r}
~,
\]
which are induced by the universal property of limits.

\item As the suspension spectrum functor is symmetric monoidal, the diagram
\begin{equation}
\label{susp.is.symm.mon}
\begin{tikzcd}[row sep=1.5cm, column sep=1.5cm]
\int_{|\cyclic} \BS^\fin
\arrow{r}{\prod}
\arrow{d}[swap]{\int_{|\cyclic} \fB \Sigma^\infty_+}
&
\Spaces
\arrow{d}{\Sigma^\infty_+}
\\
\int_{|\cyclic} \BSp^\fin
\arrow{r}[swap]{\bigotimes}
&
\Spectra
\end{tikzcd}
\end{equation}
canonically commutes.  This explains the first association
\[
\left( \prod,\htpy \right)^\fin
\longmapsto
\left( \Sigma^\infty_+ \circ \prod , \htpy \right)^\fin
\simeq
\left( \bigotimes \circ \int_{|\cyclic} \fB \Sigma^\infty_+ , \htpy \right)^\fin
\]
of \Cref{tate.package.from.diagonal.package}.

\item It is immediate that both (equivalent) composites in commutative square \Cref{susp.is.symm.mon} lie in the subcategory
\[
\Lin \left( \int_{|\cyclic} \BS^\fin , \Spectra \right)
\subset
\Fun \left( \int_{|\cyclic} \BS^\fin , \Spectra \right)
~,
\]
as the upper functor $\prod$ is fiberwise multi-cocontinuous while the right functor $\Sigma^\infty_+$ is linear.  Thus, the second association indeed factors as indicated.

\item The restriction
\begin{equation}
\label{restriction.on.fiberwise.multilinear.functors.to.Sp}
\Lin \left( \int_{|\cyclic} \BS^\fin , \Spectra \right)
\xla{- \circ \left( \int_{|\cyclic} \fB \Sigma^\infty_+\right)}
\Lin \left( \int_{|\cyclic} \BSp^\fin , \Spectra \right)
\end{equation}
is an equivalence due to a multivariate elaboration of the universal property of $\Spectra^\fin$ as the \textit{stabilization} of $\Spaces^\fin$, namely that the suspension spectrum functor
\[
\Spaces^\fin
\xra{\Sigma^\infty_+}
\Spectra^\fin
\]
is the initial linear functor from $\Spaces^\fin$ into a \textit{stable} $\infty$-category (considered in $\StCat$, i.e.\! with respect to exact functors between stable $\infty$-categories).\footnote{To see this multivariate elaboration, we work fiberwise over $\cyclic$.  Namely, writing $\cV$ simultaneously for $\Spaces^\fin$ and $\Spectra^\fin$, let us consider the full subcategory
\[
\Lin^\rel_{/\cyclic} \left( \int_{|\cyclic}\BV, \ul{\Spectra} \right)
\subset
\Fun^\rel_{/\cyclic} \left( \int_{|\cyclic}\BV, \ul{\Spectra} \right)
\]
on the multilinear functors.  The restriction \Cref{restriction.on.fiberwise.multilinear.functors.to.Sp} arises as the global sections of a morphism
\[ \begin{tikzcd}[ampersand replacement=\&, column sep=0.5cm]
\Lin^\rel_{/\cyclic} \left( \int_{|\cyclic} \BS^\fin , \ul{\Spectra} \right)
\arrow{rd}
\&
\&
\Lin^\rel_{/\cyclic} \left( \int_{|\cyclic} \BSp^\fin , \ul{\Spectra} \right)
\arrow{ll}
\arrow{ld}
\\
\&
\cyclic
\end{tikzcd} \]
in $\Cat_{/\cyclic}$, which it therefore suffices to show is an equivalence.  But this morphism evidently lies in the subcategory $\Cart_\cyclic \subset \Cat_{/\cyclic}$, and so it suffices to show that it is an equivalence on fibers.  But this is clear.}  Moreover, the forgetful functor
\[
\RMod_{\llax.\Nx}
\longra
\Cat
\]
from left-lax right $\Nx$-modules (and \textit{strict} $\Nx$-equivariant functors) to their underlying $\infty$-categories is conservative.

\end{enumerate}

Now, to obtain the Tate package from the finite Tate package, we observe that the functor \Cref{cocart.over.BW.to.llax.r.Nx.mod} takes the morphism
\begin{equation}
\label{include.BSpfin.into.BSp.over.cyclic.and.BW}
\begin{tikzcd}
\int_{|\cyclic} \BSp^\fin
\arrow[hook]{rr}{\ff}
\arrow{rd}
&
&
\int_{|\cyclic} \BSp
\arrow{ld}
\\
&
\cyclic
\arrow{d}
\\
&
\BW
\end{tikzcd}
\end{equation}
in $\coCart_\BW$ to a morphism
\[
\left(
\Fun \left(
\int_{|\cyclic} \BSp^\fin
,
\Spectra
\right)
\ractstau
\Nx
\right)
\longla
\left(
\Fun \left(
\int_{|\cyclic} \BSp
,
\Spectra
\right)
\ractstau
\Nx
\right)
\]
in $\RMod_{\llax.\Nx}$ which is a right adjoint on underlying $\infty$-categories, and so extends to an adjunction
\begin{equation}
\label{adjn.in.RMod.llax.Nx.going.from.finite.tate.package.to.tate.package}
\begin{tikzcd}[column sep=2cm]
\left(
\Fun \left(
\int_{|\cyclic} \BSp^\fin
,
\Spectra
\right)
\ractstau
\Nx
\right)
\arrow[dashed, transform canvas={yshift=0.9ex}]{r}
\arrow[leftarrow, transform canvas={yshift=-0.9ex}]{r}[transform canvas={yshift=-0.2ex}]{\bot}
&
\left(
\Fun \left(
\int_{|\cyclic} \BSp
,
\Spectra
\right)
\ractstau
\Nx
\right)
\end{tikzcd}
\end{equation}
in $\RMod^\rlax_{\llax.\Nx}$ [\Cref{cyclo:lem.get.r.lax.left.adjt}].  Taking right-lax limits of this left adjoint, we obtain a functor
\[ \begin{tikzcd}[row sep=0cm, column sep=2cm]
\lim^\rlax \left(
\Fun \left(
\int_{|\cyclic} \BSp^\fin
,
\Spectra
\right)
\ractstau
\Nx
\right)
\arrow{r}
&
\lim^\rlax \left(
\Fun \left(
\int_{|\cyclic} \BSp
,
\Spectra
\right)
\ractstau
\Nx
\right)
\\
\rotatebox{90}{$\in$}
&
\rotatebox{90}{$\in$}
\\
\left( \bigotimes , \tate \right)^\fin
\arrow[maps to]{r}
&
\left( \bigotimes , \tate \right)
\end{tikzcd} \]
taking the finite Tate package to the Tate package.

To see that this object indeed overlies the ``tensor everything together'' functor \Cref{tensor.everything.together.over.cyclic}, we first observe that on underlying $\infty$-categories, the adjunction \Cref{adjn.in.RMod.llax.Nx.going.from.finite.tate.package.to.tate.package} in $\RMod^\rlax_{\llax.\Nx}$ recovers an adjunction
\[ \begin{tikzcd}[column sep=2cm]
\Fun \left(
\int_{|\cyclic} \BSp^\fin
,
\Spectra
\right)
\arrow[dashed, transform canvas={yshift=0.9ex}]{r}
\arrow[leftarrow, transform canvas={yshift=-0.9ex}]{r}[transform canvas={yshift=-0.2ex}]{\bot}
&
\Fun \left(
\int_{|\cyclic} \BSp
,
\Spectra
\right)
\end{tikzcd} \]
in $\Cat$.  This left adjoint is given by left Kan extension along the functor \Cref{include.BSpfin.into.BSp.over.cyclic.and.BW} (i.e.\! along the image in $\Cat$ of that morphism in $\coCart_\BW$); considering that functor as a morphism in $\coCart_\cyclic$, we see that this left Kan extension may be computed fiberwise over $\cyclic$.  We then conclude by noting that for any $r \geq 1$, the commutative triangle
\[ \begin{tikzcd}[column sep=1.5cm]
\left( \Spectra^\fin \right)^{\times r}
\arrow{r}{\bigotimes}
\arrow[hook]{d}[swap]{\ff}
&
\Spectra
\\
\Spectra^{\times r}
\arrow{ru}[swap]{\bigotimes}
\end{tikzcd} \]
is a left Kan extension diagram: for any $E \in \Spectra$ the $\infty$-category
\[
\left( \Spectra^\fin \right)_{/E}
:=
\Spectra^\fin
\underset{\Spectra}{\times}
\Spectra_{/E}
\]
is sifted (in fact filtered), so for any $(E_i) \in \Spectra^{\times r}$ the $\infty$-category
\[
\left( \left( \Spectra^\fin \right)^{\times r} \right)_{/(E_i)}
:=
\left( \Spectra^\fin \right)^{\times r}
\underset{\Spectra^{\times r}}{\times}
\left( \Spectra^{\times r} \right)_{/(E_i)}
\simeq
\prod_i
\left(
\Spectra^\fin
\underset{\Spectra}{\times}
\Spectra_{/E_i}
\right)
\]
is sifted as well, and the tensor product functor
\[
\Spectra^{\times r}
\xlongra{\bigotimes}
\Spectra
\]
commutes with sifted colimits.

\subsection{The cyclotomic structure on $\THH$}
\label{subsection.cyclo.str.on.THH}

By \Cref{cyclo.is.rlax.lim}, we can identify the $\infty$-category of cyclotomic spectra as the right-lax limit
\[
\Cyclo(\Spectra)
\simeq
\lim^\rlax \left( 
\Fun( \BT , \Spectra )
\ractstau
\Nx
\right)
~.
\]
Then, for any spectrally-enriched $\infty$-category $\cC$, the Tate package is taken to $\THH(\cC)$ equipped with its cyclotomic structure as illustrated in \Cref{Tate.package.gives.cyclo.str}.
\begin{sidewaysfigure}
\vspace{425pt}
\[ \begin{tikzcd}[row sep=2.5cm, column sep=2cm]
\Fun \left( \int_{|\cyclic} \BSp, \Spectra \right)
\arrow{r}
\arrow[leftarrow]{d}
&
\Fun \left( \int_{|\cyclic} \cdiC , \Spectra \right)
\arrow[dashed, transform canvas={yshift=0.9ex}]{r}
\arrow[leftarrow, transform canvas={yshift=-0.9ex}]{r}[transform canvas={yshift=-0.2ex}]{\bot}
\arrow[leftarrow]{d}
&
\Fun ( \BT , \Spectra )
\arrow[leftarrow]{d}
\\
\lim^\rlax
\left(
\Fun \left( \int_{|\cyclic} \BSp , \Spectra \right)
\ractstau
\Nx
\right)
\arrow{r}
&
\lim^\rlax
\left(
\Fun \left( \int_{|\cyclic} \cdiC , \Spectra \right)
\ractstau
\Nx
\right)
\arrow[dashed, transform canvas={yshift=0.9ex}]{r}
\arrow[leftarrow, transform canvas={yshift=-0.9ex}]{r}[transform canvas={yshift=-0.2ex}]{\bot}
&
\lim^\rlax
\left(
\Fun \left( \BT , \Spectra \right)
\ractstau
\Nx
\right)
&
\hspace{-2.6cm} \simeq \Cyclo(\Spectra)
\\
\bigotimes
\arrow[mapsto]{rr}{\Cref{fact.hlgy.with.T.action}}
&
&
\left( \TT \lacts \THH(\cC) \right)
\\
\left( \bigotimes , \tate \right)
\arrow[mapsto]{u}
\arrow[mapsto]{rr}
&
&
\left( \Cyclo \lacts \THH(\cC) \right)
\arrow[mapsto]{u}
\end{tikzcd} \]
\vspace{2.5cm}
\caption{Above: the diagram of $\infty$-categories housing the passage from the Tate package to the cyclotomic structure on $\THH(\cC)$.  All solid horizontal morphisms are given by pullback; the upwards functors are given by restriction to underlying objects of left-lax right $\Nx$-modules, and the diagram commutes with omitting either both left adjoints or both right adjoints.  Below: the Tate package gives rise to the cyclotomic structure on $\THH(\cC)$, and the restriction of this passage recovers its underlying $\TT$-action.}
\label{Tate.package.gives.cyclo.str}
\end{sidewaysfigure}
In analogy with \Cref{unstable.cyclo.str.via.modules}, we also summarize the construction in \Cref{cyclo.str.via.modules}.
\begin{figure}[h]
\begin{tikzcd}[column sep=4.5cm, row sep=1.5cm]
{[1] \underset{1,\pt,1}{\coprod} [1]}
\arrow{r}{\left( \int_{|\epicyc} \BSp \la \int_{|\epicyc} \cdiC \ra \BW \right)^\circ}
\arrow{d}[swap]{\id \underset{\pt}{\coprod} {\sf r.adjt}}  
&
\left( \LMod_\WW \right)^\op
\arrow{r}{\left( \Fun ( (-)_0 , \Spectra) \ractstau \Nx \right)}
&
\RMod_{\llax.\Nx}
\arrow[hook, two heads]{d}
\\
{ [1] \underset{1,\pt,l}{\coprod} \Adj}
\arrow[dashed]{rr}{\exists ! ~ [\textup{\Cref{cyclo:lem.get.r.lax.left.adjt}}]}
\arrow[leftarrow]{d}[swap]{\id \underset{\pt}{\coprod} {\sf l.adjt}}  
&
&
\RMod^\rlax_{\llax.\Nx}
\arrow{dd}{\lim^\rlax}
\\
{[2]}
\arrow[leftarrow]{d}[swap]{\{0 < 2 \}}   
\\
{[1]}
\arrow[dashed]{rr}
&
&
\Cat
\end{tikzcd}
\caption{A diagram illustrating in the language of modules how to obtain the lower rightwards composite functor in \Cref{Tate.package.gives.cyclo.str}, which takes the Tate package to the cyclotomic structure on $\THH(\cC)$.  Here, we write $l \in \Adj$ for the source of the walking left adjoint.}
\label{cyclo.str.via.modules}
\end{figure}

Unlike in \Cref{diagonal.package.gives.W.action}, the commutativity of the upper diagram in \Cref{Tate.package.gives.cyclo.str} (after omitting either both left adjoints or both right adjoints) -- which guarantees that the underlying spectrum of the cyclotomic spectrum that is constructed is indeed $\THH(\cC)$ -- is not trivial.  However, it can be seen as follows.  From [\Cref{cyclo:define.rest.of.equivariant.functors} and \Cref{cyclo:obs.could.use.Delta.over.B.or.Cat.over.B}], we have a (nearly definitionally) commutative diagram
\[ \begin{tikzcd}[row sep=1.5cm]
\RMod^\rlax_{\llax.\Nx}
\arrow{rr}{\lim^\rlax}
\arrow[hook]{rd}[sloped, swap, pos=0.6]{\ff}
&
&
\Cat
\\
&
\Fun \left( \left( \Cat_{/\BNop} \right)^\op , \Cat \right)
\arrow{ru}[swap]{\ev_{\sd(\BNop)}}
\end{tikzcd}~. \]
Note that evaluating the image of a left-lax right $\Nx$-module on a morphism
\begin{equation}
\label{morphism.in.Cat.over.BN.op.giving.map.from.rlax.lim.to.underlying.obj}
\begin{tikzcd}
\pt
\arrow{rr}
\arrow{rd}
&
&
\sd(\BNop)
\arrow{ld}
\\
&
\BNop
\end{tikzcd}
\end{equation}
in $\Cat_{/\BNop}$ determines the canonical map from its right-lax limit to its underlying object -- the upwards functors in \Cref{Tate.package.gives.cyclo.str}.  Then, the commutativity of the left square in \Cref{Tate.package.gives.cyclo.str} is obtained by precomposing the adjoint of the composite
\[
[1]
\xra{\left( \int_{|\epicyc} \BSp \la \int_{|\epicyc} \cdiC \right)^\circ}
\left( \LMod_\WW \right)^\op
\xra{\left( \Fun ( (-)_0 , \Spectra ) \ractstau \Nx \right)}
\RMod_{\llax.\Nx}
\longhookra
\Fun \left( \left( \Cat_{/\BNop} \right)^\op , \Cat \right)
\]
to obtain the composite
\[
[1] \times [1]
\xra{\id \times \Cref{morphism.in.Cat.over.BN.op.giving.map.from.rlax.lim.to.underlying.obj}^\circ}
[1] \times \left( \Cat_{/\BNop} \right)^\op
\longra
\Cat
~.
\]
Similarly, the commutativity of the right square in \Cref{Tate.package.gives.cyclo.str} is obtained by precomposing the adjoint of the composite
\[
\Adj
\xra{
\scalebox{0.8}{
\begin{tikzcd}[column sep=1.5cm, ampersand replacement=\&]
\left( \Fun \left( \int_{|\cyclic} \cdiC , \Spectra \right) \ractstau \Nx \right)
\arrow[transform canvas={yshift=0.9ex}]{r}
\arrow[leftarrow, transform canvas={yshift=-0.9ex}]{r}[transform canvas={yshift=-0.2ex}]{\bot}
\&
\left( \Fun ( \BT , \Spectra) \ractstau \Nx \right)
\end{tikzcd}
}
}
\RMod^\rlax_{\llax.\Nx}
\overset{\ff}{\longhookra}
\Fun \left( \left( \Cat_{/\BNop} \right)^\op , \Cat \right)
\]
to obtain the composite
\[
\Adj \times [1]
\xra{\id \times \Cref{morphism.in.Cat.over.BN.op.giving.map.from.rlax.lim.to.underlying.obj}^\circ}
\Adj \times \left( \Cat_{/\BNop} \right)^\op
\longra
\Cat
~.
\eqno \qed
\]

\part{The cyclotomic trace}
\label{part.cyclo.trace}

In this part of the paper, we prove \Cref{cyclo.trace} modulo certain key maneuvers which are contained in the papers \cite{AMR-cyclo,AMR-fact}.  To begin, in \Cref{section.unstable.cyclo.trace} we recall the unstable cyclotomic trace for flagged spatially-enriched $\infty$-categories, as constructed in [\Cref{fact:section.unstable.cyclo.trace}].  Then, in \Cref{section.cyclo.trace}, through a key maneuver which ultimately relies on the fundamental calculus-theoretic feature of the Tate construction [\Cref{cyclo:Cr.genzd.tate.of.rth.power.is.exact}], we use this to obtain a cyclotomic pre-trace map for spectrally-enriched $\infty$-categories and thereafter the cyclotomic trace map for stable $\infty$-categories.

\section{The unstable cyclotomic trace}
\label{section.unstable.cyclo.trace}

Factorization homology enriched in a cartesian symmetric monoidal $\infty$-category $(\cV,\times)$ assembles into a bifunctor
\[
\int_{(-)} (-)
:
\M \times \fCat(\cV)
\longra
\cV
\]
[\Cref{fact:definition.fully.functorial.enr.fact.hlgy} \and \Cref{fact:D.of.to.D.over.is.final}].  Note the fully faithful embedding
\[
\BW^\lcone
\overset{\ff}{\longhookra}
\M
\]
extending the defining inclusion $\BW \hookra \M$ over the initial object $\DD^0 \in \M$.  This determines the \bit{unstable cyclotomic trace} map
\begin{equation}
\label{unstable.cyclo.trace.as.functors}
\int_{\DD^0}
\simeq
\iota
\longra
\TCSx
:=
\left( \THHVx \right)^{\htpy \WW}
~,
\end{equation}
which assembles into a natural transformation
\[ \begin{tikzcd}[column sep=1.5cm]
\fCat(\Spaces)
\arrow[bend left]{r}[pos=0.4]{\iota}[swap, transform canvas={yshift=-0.25cm}]{\Downarrow}   
\arrow[bend right]{r}[swap, pos=0.4]{\TCSx}
&
\Spaces
\end{tikzcd} ~. \]

\section{The cyclotomic trace}
\label{section.cyclo.trace}

Given any flagged spectrally-enriched $\infty$-category $\cC \in \fCat(\Spectra)$, consider its underlying $\infty$-groupoid
\[
\iC
\in
\Spaces
\subset
\Cat
\simeq
\Cat(\Spaces)
\subset
\fCat(\Spaces)
\]
as a flagged spatially-enriched $\infty$-category, with enriched hom functor
\[
\cdiC
\xra{\ulhom_\iC}
\BS
~.
\]
The postcomposition
\[
\ulhom_\BSiC
:
\cdiC
\xra{\ulhom_\iC}
\BS
\xra{\fB\Sigma^\infty_+}
\BSp
\]
then determines the enriched hom functor of a flagged spectrally-enriched $\infty$-category
\[
\BSiC \in \fCat(\Spectra)
~,
\]
which is the initial flagged spectrally-enriched $\infty$-category whose underlying $\infty$-groupoid is $\iC$.\footnote{To see this, note first that $\iC$ is evidently the initial flagged spatially-enriched $\infty$-category with underlying $\infty$-groupoid given by $\iC$ (e.g.\! using the equivalence between flagged spatially-enriched $\infty$-categories and Segal spaces \cite[Theorem 4.4.7]{GH-enr}).  Then, observe that in the adjunction $\Sigma^\infty_+ \adj \Omega^\infty$ the left adjoint is symmetric monoidal while the right adjoint is left-lax symmetric monoidal, so that taking right-lax functors of category objects out of $\cdiC$ gives another adjunction, whose left adjoint (which preserves initial objects) takes $\iC$ to $\BSiC$.}  Hence, we obtain a canonical functor
\begin{equation}
\label{canonical.morphism.from.BSiC.to.C}
\BSiC
\longra
\cC
\end{equation}
in $\fCat(\Spectra)$.  This determines the \textit{cyclotomic pre-trace}, namely the composite morphism
\begin{align}
\nonumber
\Sigma^\infty_+(\iC)
&
\xra{\Sigma^\infty_+ \Cref{unstable.cyclo.trace.as.functors}}
\Sigma^\infty_+(\TCSx(\iC))
\\
\nonumber
&
=:
\Sigma^\infty_+ \left( \THHSx(\iC)^{\htpy \WW} \right)
\\
\label{commute.Sigma.infty.plus.with.htpy.W.fixedpoints}
&
\longra
\left( \Sigma^\infty_+ \THHSx(\iC) \right)^{\htpy \WW}
\\
\label{canonical.map.for.THHSx.iC.as.a.cyclo.spectrum.with.Frob.lifts}
&
\longra
\left( \Sigma^\infty_+ \THHSx(\iC) \right)^{\htpy \Cyclo}
\\
\label{equivalence.between.THHSxiC.and.THHSpBSiC}
&
\simeq
\left( \THH_\Spectra ( \BSiC ) \right)^{\htpy \Cyclo}
\\
\nonumber
&
\xra{\left( \THH_\Spectra \Cref{canonical.morphism.from.BSiC.to.C} \right)^{\htpy \Cyclo}}
\left( \THH_\Spectra(\cC) \right)^{\htpy \Cyclo}
\\
\nonumber
&
=:
\TC(\cC)
\end{align}
in $\Spectra$, as we now explain.\footnote{A similar maneuver (or just postcomposition with the canonical map $\TC(\cC) \ra \THH(\cC)$) yields the \textit{Dennis pre-trace}, a morphism $\Sigma^\infty_+(\iC) \ra \THH(\cC)$.}
\begin{itemize}

\item The morphism \Cref{commute.Sigma.infty.plus.with.htpy.W.fixedpoints} arises from the universal property of limits.

\item To obtain the morphism \Cref{canonical.map.for.THHSx.iC.as.a.cyclo.spectrum.with.Frob.lifts}, recall first that $\THHSx(\iC)$ is an unstable cyclotomic space, so that $\Sigma^\infty_+ \THHSx(\iC)$ is an unstable cyclotomic \textit{spectrum}: the composite
\[
\BW
\xra{\THHSx(\iC)}
\Spaces
\xra{\Sigma^\infty_+}
\Spectra
~.
\]
By [\Cref{cyclo:cor.cart.fibn.to.BN.for.htpy.cyclo.sp}], this may also be identified as a \bit{cyclotomic spectrum with Frobenius lifts} [\Cref{cyclo:define.cyclo.spt.with.frob.lifts}]: a cyclotomic spectrum $T \in \Cyclo(\Spectra)$ equipped with suitably equivariant lifts
\[ \begin{tikzcd}[row sep=1.5cm, column sep=1.5cm]
T
\arrow[dashed]{r}{\tilde{\sigma}_r}
\arrow{rd}[swap]{\sigma_r}
&
T^{\htpy \Cyclic_r}
\arrow{r}
\arrow{d}
&
T
\\
&
T^{\tate \Cyclic_r}
\end{tikzcd}~, \]
compatibly for all $r \in \Nx$.  Then, the morphism \Cref{canonical.map.for.THHSx.iC.as.a.cyclo.spectrum.with.Frob.lifts} is guaranteed by [\Cref{cyclo:obs.h.to.tau.map.for.cyclo.spt.w.frob.lifts}].

\item The equivalence \Cref{equivalence.between.THHSxiC.and.THHSpBSiC} results from the fact that the suspension spectrum functor
\[
\Spaces
\xra{\Sigma^\infty_+}
\Spectra
\]
is symmetric monoidal and commutes with colimits.

\end{itemize}
This construction is natural in the variable $\cC \in \fCat(\cC)$, so that it determines the natural transformation in the diagram
\[ \begin{tikzcd}[column sep=1.5cm]
\StCat
\arrow{r}
&
\fCat(\Spectra)
\arrow[bend left]{r}[pos=0.45]{\Sigma^\infty_+ \iota}[swap, transform canvas={yshift=-0.3cm}]{\rotatebox{-90}{$\Rightarrow$}}
\arrow[bend right]{r}[swap, pos=0.45]{\TC}
&
\Spectra
\end{tikzcd}~. \]
By the universal property of algebraic K-theory for stable $\infty$-categories \cite{BGT-K}, this induces a unique factorization
\[ \begin{tikzcd}[column sep=4.5cm]
\StCat
\arrow[bend left=50]{r}{\Sigma^\infty_+ \iota}
\arrow{r}[pos=0.47]{\K}[transform canvas={yshift={0.55cm}}, pos=0.47]{\Downarrow}[swap, transform canvas={yshift={-0.55cm}}, pos=0.47]{\Downarrow}
\arrow[bend right=50]{r}[swap]{\TC}
&
\Spectra
\end{tikzcd}~, \]
namely the cyclotomic trace. \hfill \qed

\addtocontents{toc}{\protect\vspace{10pt}}

\bibliographystyle{amsalpha}
\bibliography{trace}{}

\providecommand{\bysame}{\leavevmode\hbox to3em{\hrulefill}\thinspace}
\providecommand{\MR}{\relax\ifhmode\unskip\space\fi MR }
\providecommand{\MRhref}[2]{%
  \href{http://www.ams.org/mathscinet-getitem?mr=#1}{#2}
}
\providecommand{\href}[2]{#2}
\begin{thebibliography}{AGHL14}

\bibitem[AF]{AF-fibns}
David Ayala and John Francis, \emph{Fibrations of $\infty$-categories},
  available at \texttt{arXiv:1702.02681}, v1.

\bibitem[AFR]{AFR-fact}
David Ayala, John Francis, and Nick Rozenblyum, \emph{Factorization homology
  {I}: higher categories}, available at \texttt{arXiv:1504.04007}, v4.

\bibitem[AFT17]{AFT-fact}
David Ayala, John Francis, and Hiro~Lee Tanaka, \emph{Factorization homology of
  stratified spaces}, Selecta Math. (N.S.) \textbf{23} (2017), no.~1, 293--362.

\bibitem[AGH09]{AngGerHess-trunc}
Vigleik Angeltveit, Teena Gerhardt, and Lars Hesselholt, \emph{On the
  {$K$}-theory of truncated polynomial algebras over the integers}, J. Topol.
  \textbf{2} (2009), no.~2, 277--294.

\bibitem[AGHL14]{AngGerHillLind-trunc}
Vigleik Angeltveit, Teena Gerhardt, Michael~A. Hill, and Ayelet Lindenstrauss,
  \emph{On the algebraic {$K$}-theory of truncated polynomial algebras in
  several variables}, J. K-Theory \textbf{13} (2014), no.~1, 57--81.

\bibitem[AMGRa]{AMR-fact}
David Ayala, Aaron Mazel-Gee, and Nick Rozenblyum, \emph{Factorization homology
  of enriched $\infty$-categories}, appearing concurrently.

\bibitem[AMGRb]{AMR-cyclo}
\bysame, \emph{A naive approach to genuine {$G$}-spectra and cyclotomic
  spectra}, appearing concurrently.

\bibitem[BCD10]{BCD-cov}
Morten Brun, Gunnar Carlsson, and Bj\o rn~Ian Dundas, \emph{Covering homology},
  Adv. Math. \textbf{225} (2010), no.~6, 3166--3213.

\bibitem[BGT13]{BGT-K}
Andrew~J. Blumberg, David Gepner, and Gon\c{c}alo Tabuada, \emph{A universal
  characterization of higher algebraic {$K$}-theory}, Geom. Topol. \textbf{17}
  (2013), no.~2, 733--838.

\bibitem[BHM93]{BHM}
M.~B{\"o}kstedt, W.~C. Hsiang, and I.~Madsen, \emph{The cyclotomic trace and
  algebraic {$K$}-theory of spaces}, Invent. Math. \textbf{111} (1993), no.~3,
  465--539.

\bibitem[BM]{BluMan-AKTS}
Andre Blumberg and Michael Mandell, \emph{The homotopy groups of the algebraic
  {K}-theory of the sphere spectru,}, available at \texttt{arXiv:1408.0133},
  v2.

\bibitem[BM12]{BluMan-loc}
Andrew~J. Blumberg and Michael~A. Mandell, \emph{Localization theorems in
  topological {H}ochschild homology and topological cyclic homology}, Geom.
  Topol. \textbf{16} (2012), no.~2, 1053--1120.

\bibitem[BS58]{BS-GRR}
Armand Borel and Jean-Pierre Serre, \emph{Le th\'eor\`eme de {R}iemann-{R}och},
  Bull. Soc. Math. France \textbf{86} (1958), 97--136.

\bibitem[BZFN10]{BZFN}
David Ben-Zvi, John Francis, and David Nadler, \emph{Integral transforms and
  {D}rinfeld centers in derived algebraic geometry}, J. Amer. Math. Soc.
  \textbf{23} (2010), no.~4, 909--966.

\bibitem[Con72]{Cond-fibns}
Fran\c{c}ois Conduch\'e, \emph{Au sujet de l'existence d'adjoints \`a droite
  aux foncteurs ``image r\'eciproque'' dans la cat\'egorie des cat\'egories},
  C. R. Acad. Sci. Paris S\'er. A-B \textbf{275} (1972), A891--A894.

\bibitem[DGM13]{DGM-book}
Bj\o rn~Ian Dundas, Thomas~G. Goodwillie, and Randy McCarthy, \emph{The local
  structure of algebraic {K}-theory}, Algebra and Applications, vol.~18,
  Springer-Verlag London, Ltd., London, 2013.

\bibitem[DM94]{DM-KTHH}
Bj\o rn~Ian Dundas and Randy McCarthy, \emph{Stable {$K$}-theory and
  topological {H}ochschild homology}, Ann. of Math. (2) \textbf{140} (1994),
  no.~3, 685--701.

\bibitem[Dun97]{DundasRelAKT}
Bj{\o}rn~Ian Dundas, \emph{Relative {$K$}-theory and topological cyclic
  homology}, Acta Math. \textbf{179} (1997), no.~2, 223--242.

\bibitem[FG05]{K-handbook}
Eric~M. Friedlander and Daniel~R. Grayson (eds.), \emph{Handbook of
  {$K$}-theory. {V}ol. 1, 2}, Springer-Verlag, Berlin, 2005.

\bibitem[GH15]{GH-enr}
David Gepner and Rune Haugseng, \emph{Enriched {$\infty$}-categories via
  non-symmetric {$\infty$}-operads}, Adv. Math. \textbf{279} (2015), 575--716.

\bibitem[Gir64]{Gir-desc}
Jean Giraud, \emph{M\'ethode de la descente}, Bull. Soc. Math. France M\'em.
  \textbf{2} (1964), viii+150.

\bibitem[Goo86]{GooRelAKT}
Thomas~G. Goodwillie, \emph{Relative algebraic {$K$}-theory and cyclic
  homology}, Ann. of Math. (2) \textbf{124} (1986), no.~2, 347--402.

\bibitem[GR17]{GR}
Dennis Gaitsgory and Nick Rozenblyum, \emph{{A Study in Derived Algebraic
  Geometry, Volume I: Correspondences and Duality}}, Mathematical Surveys and
  Monographs, vol. 221, American Mathematical Society, Providence, RI, 2017.

\bibitem[HM97a]{HessMad-trunc}
Lars Hesselholt and Ib~Madsen, \emph{Cyclic polytopes and the {$K$}-theory of
  truncated polynomial algebras}, Invent. Math. \textbf{130} (1997), no.~1,
  73--97.

\bibitem[HM97b]{HessMad-Witt}
\bysame, \emph{On the {$K$}-theory of finite algebras over {W}itt vectors of
  perfect fields}, Topology \textbf{36} (1997), no.~1, 29--101.

\bibitem[HM03]{HessMad-local}
\bysame, \emph{On the {$K$}-theory of local fields}, Ann. of Math. (2)
  \textbf{158} (2003), no.~1, 1--113.

\bibitem[HM04]{HessMad-dRW}
\bysame, \emph{On the {D}e {R}ham-{W}itt complex in mixed characteristic}, Ann.
  Sci. \'Ecole Norm. Sup. (4) \textbf{37} (2004), no.~1, 1--43.

\bibitem[KR97]{KleinRog-fib}
John~R. Klein and John Rognes, \emph{The fiber of the linearization map
  {$A(*)\to K({\bf Z})$}}, Topology \textbf{36} (1997), no.~4, 829--848.

\bibitem[Lod98]{Loday-cyclic}
Jean-Louis Loday, \emph{Cyclic homology}, second ed., Grundlehren der
  Mathematischen Wissenschaften [Fundamental Principles of Mathematical
  Sciences], vol. 301, Springer-Verlag, Berlin, 1998, Appendix E by Mar\'\i a
  O. Ronco, Chapter 13 by the author in collaboration with Teimuraz Pirashvili.

\bibitem[Lur]{LurieHA}
Jacob Lurie, \emph{Higher {A}lgebra}, available from the author's website
  (version dated May 16, 2016).

\bibitem[Lur09]{LurieHTT}
\bysame, \emph{Higher topos theory}, Annals of Mathematics Studies, vol. 170,
  Princeton University Press, Princeton, NJ, 2009.

\bibitem[McC97]{McCRelAKT}
Randy McCarthy, \emph{Relative algebraic {$K$}-theory and topological cyclic
  homology}, Acta Math. \textbf{179} (1997), no.~2, 197--222.

\bibitem[MSV97]{MSV-THH}
J.~McClure, R.~Schw{\"a}nzl, and R.~Vogt, \emph{{$THH(R)\cong R\otimes S^1$}
  for {$E_\infty$} ring spectra}, J. Pure Appl. Algebra \textbf{121} (1997),
  no.~2, 137--159.

\bibitem[NS]{NS}
Thomas Nikolaus and Peter Scholze, \emph{On topological cyclic homology},
  available at \texttt{arXiv:1707.01799}, v1.

\bibitem[Rog03]{Rog-White}
John Rognes, \emph{The smooth {W}hitehead spectrum of a point at odd regular
  primes}, Geom. Topol. \textbf{7} (2003), 155--184.

\bibitem[SGA71]{SGA6}
\emph{Th\'eorie des intersections et th\'eor\`eme de {R}iemann-{R}och}, Lecture
  Notes in Mathematics, Vol. 225, Springer-Verlag, Berlin-New York, 1971,
  S\'eminaire de G\'eom\'etrie Alg\'ebrique du Bois-Marie 1966--1967 (SGA 6),
  Dirig\'e par P. Berthelot, A. Grothendieck et L. Illusie. Avec la
  collaboration de D. Ferrand, J. P. Jouanolou, O. Jussila, S. Kleiman, M.
  Raynaud et J. P. Serre.

\bibitem[TV09]{TV-loops}
Bertrand To{\"e}n and Gabriele Vezzosi, \emph{Chern character, loop spaces and
  derived algebraic geometry}, Algebraic topology, Abel Symp., vol.~4,
  Springer, Berlin, 2009, pp.~331--354.

\end{thebibliography}

\end{document}